
%
%
\documentstyle{amsppt}
\pagewidth{6.4in}\vsize8.5in\parindent=6mm
\parskip=3pt\baselineskip=14pt\tolerance=10000\hbadness=500
\NoRunningHeads
\loadbold
\topmatter
\title
Sharp Lorentz space estimates for rough operators
\endtitle
\author
Andreas Seeger and Terence Tao
\endauthor
\abstract
We demonstrate the $(H^1,L^{1,2})$ or $(L^p,L^{p,2})$ mapping properties of several rough operators.  In all cases these estimates are sharp in the sense that the Lorentz exponent $2$ cannot be replaced by any lower number.
\endabstract
\thanks 
This work was conducted at Madison, Princeton, UCLA, and UNSW.
The research was  supported in part by  grants from the National Science 
Foundation (DMS 9501040, DMS 9970042 and DMS 9706764).
  \endthanks
\address 
Andreas Seeger, University of Wisconsin, Madison, WI 53706-1388
\endaddress 
\email seeger\@math.wisc.edu\endemail
\address 
Terence Tao,
University of California, Los Angeles, CA  90095-1555, 
and
 School of Mathematics, University of New South Wales, 
Sydney NSW 2052, Australia
\endaddress 
\email tao\@math.ucla.edu\endemail 
\endtopmatter
\document


\def\R{\Bbb R}

\def\diam{{\text{\rm diam}}}

\define\ffB{\Cal B}
\def\itemize#1{\item"{#1}"}
\def\seq{\subseteq}
\define\Id{\text{\sl Id}}

\define\Ga{\Gamma}
\define\ga{\gamma}
\redefine\Re{\text{\rm Re}}
\redefine\Im{\text{\rm Im}}

\define\prd{{\text{\it prod}}}
\define\parab{{\text{\it parabolic}}}
\define\vth{\vartheta}


\define\eg{{\it e.g. }}

\define\sgn{{\text{\rm sign }}}


\define\dist{{\text{\rm dist}}}

\define\inn#1#2{\langle#1,#2\rangle}

\define\lcontr{\rfloor}
\define\lco#1#2{{#1}\lcontr{#2}}
\define\lcoi#1#2{\imath({#1}){#2}}
\define\rco#1#2{{#1}\rcontr{#2}}
\redefine\exp{{\text{\rm exp}}}
\define\bin#1#2{{\pmatrix {#1}\\{#2}\endpmatrix}}
\define\meas{{\text{\rm meas}}}

\define\card{\text{\rm card}}
\define\lc{\lesssim}
\define\gc{\gtrsim}
\define\pv{\text{\rm p.v.}}


\define\del{\delta}             
\define\eps{\varepsilon}

\define\la{\lambda}             \define\La{\Lambda}

\define\fS{{\frak S}}


\define\bbR{{\Bbb R}}

\define\bbZ{{\Bbb Z}}

\define\cA{{\Cal A}}
\define\cB{{\Cal B}}

\define\cD{{\Cal D}}
\define\cE{{\Cal E}}

\define\cG{{\Cal G}}

\define\cM{{\Cal M}}

\define\cP{{\Cal P}}

\define\cR{{\Cal R}}

\define\cT{{\Cal T}}

\define\cV{{\Cal V}}
\define\cW{{\Cal W}}





\head{\bf 1. Introduction}\endhead
In this paper we consider the endpoint behaviour on Hardy spaces
of two classes of operators, 
namely singular integral operators with rough homogeneous kernels
\cite{4} and singular integral operators with  convolution kernels 
supported on curves in the plane
(\cite{20}, \cite{27}). 
These operators fall outside the Calder\'on-Zygmund theory; however
weak type  $(L^1,L^{1,\infty})$ or $(H^1,L^{1,\infty})$ inequalities  
have been established in the previous 
literature (\cite{7}, \cite{9}, \cite{16}
 \cite{18}, \cite{25}, \cite{29}) 
We shall show that 
 the target space $L^{1,\infty}$ can be improved to the 
Lorentz space $L^{1,2}$, possibly at the 
cost of moving to a stronger type of Hardy space (\eg product $H^1$).  
Examples of Christ \cite{8}, 
\cite{17}
 show that these types of results are optimal in the sense 
that one cannot replace $L^{1,2}$ 
by $L^{1,q}$ for any $q < 2$.

The space $L^{1,2}$ arises naturally as the interpolation space halfway 
between $L^{1,\infty}$ and $L^1$.  As a gross caricature of how this 
space arises, suppose that we have a collection of functions $f_i$ 
which are uniformly bounded in $L^1$, and whose maximal function 
$\sup_i |f_i|$ is in weak $L^1$, and we wish to estimate the quantity
$$ \Big\| \sum_i \gamma_i f_i \Big\|_{L^{1,2}}
$$
for some $l^2$ co-efficients $\gamma_i$.  If the $f_i$ are sufficiently 
orthogonal, we may hope to control this quantity by the square function
$$ \Big\| \big(\sum_i |\gamma_i f_i|^2\big)^{1/2} \Big\|_{L^{1,2}}.
\tag 1.1
$$
However from our hypotheses we see that
$$ \Big\| \big(\sum_i |\gamma_i f_i|^q\big)^{1/q}\Big \|_{L^{1,q}} \lesssim 
\Big(\sum_i |\gamma_i|^q\Big)^{1/q}$$
for $q=1$ and $q=\infty$, and thus by interpolation 
 for all 
$1 \leq q \leq \infty$ ({\it cf.} Lemma 2.2. below). 
 Thus we expect to control 
(1.1)
by the $\ell^2$ norm of $\{\gamma_i\}$.

Our arguments  will be based on more complicated versions 
of the above informal strategy.  Generally, the $L^1$ estimates will 
be quite trivial, whereas the $L^{1,\infty}$ estimates will be variants 
of existing weak-type (1,1)  estimates for rough operators in the 
literature (\eg \cite{7}, \cite{25}).
We shall  demonstrate this technique for two classes of 
operators.  Firstly we show that the Hilbert transform 
on plane curves $(t,t^m)$ maps product $H^1$ into $L^{1,2}$ or a related 
Hardy-Lorentz space; we also prove sharp $L^p\to L^{p,2}$ estimates for
 a related analytic family of hypersingular operators.
Then we discuss homogeneous singular integrals with rough
 kernels in $\R^d$, satisfying an $L \log^2\! L$
condition on the sphere, and show that these map the standard
 Hardy space $H^1$ to $L^{1,2}$.

We remark that a  simple  version  of the above technique has 
been used by one of the authors in \cite{23} 
to prove an endpoint version of the H\"ormander 
multiplier theorem. Namely (stating only the 
one-dimensional version) 
if $\phi$ 
is a nonzero even  smooth bump function then the condition
$\sup_{t>0}\|\phi m(t\cdot)\|_{B^2_{1/2,1}}$ implies that the 
convolution operator with Fourier multiplier 
$m$ maps $H^1$ to $L^{1,2}$  
(and an example by Baernstein and Sawyer \cite{1}
shows that $L^{1,2}$ cannot be replaced by $L^{1,q}$ for $q<2$).
The second author and Jim Wright \cite{30} 
have recently improved this result by replacing 
the Besov space
$B^2_{1/2,1}$ by the larger space  $R^2_{1/2,2}$ 
defined in \cite{24}
improving on the known  $(H^1, L^{1,\infty})$ result which is 
implicit in the latter paper.

The paper is structured as follows. After formulating
our results in the current section  we review some material
about Hardy-Lorentz spaces and interpolation, in \S2.
In \S3 we prove  an abstract variant of a  stopping time argument due to M. Christ
which may be helpful elsewhere.
\S4 contains the main square-function estimate 
needed to prove our theorems on integrals along curves;
in \S5 we conclude the proof of these results.
Rough homogeneous kernels are considered in
\S6 and \S7.

\subheading{Rough homogeneous convolution kernels}

Let $K$ be a
convolution kernel  on the Euclidean space $\R^d$ 
and assume that $K$ is homogeneous of degree $-d$ and that 
 the restriction $\Omega$ to the unit sphere is integrable and 
has mean zero, $\int_{S^{d-1}}\Omega(\theta) d\sigma(\theta)=0.$
We may  define the operator $T_\Omega$ of convolution with $K$
 on test
functions at least by the usual method of principal values:
$$
T_\Omega f(x)=p.v.\int \frac{\Omega(y/|y|)}{|y|^d} f(x-y) dy.
\tag 1.2
$$

We consider the  mapping properties of $T_\Omega$, especially near the 
endpoint $L^1$.
If $\Omega$ is somewhat regular (for example, if it is H\"older continuous
or satisfies an appropriate $L^1$ Dini condition) then
the standard Calder\'on-Zygmund theory shows that $T$ is bounded on all
$L^p$ spaces, $1 < p < \infty$, is of weak type $(1,1)$, and maps
the Hardy space $H^1$ to $L^1$.  If no regularity is assumed, but
$K$ is $L \log L$ on the sphere, then it was shown by Calder\'on-Zygmund 
\cite{4}
that $T_\Omega$ is bounded on $L^p$; in fact 
(see \cite{25}) it is of
weak type $(1,1)$.
The behaviour at $H^1$ is more subtle, however, as an example of M. Christ 
shows (see also \cite{17}).  
For the sake of
illustration let us consider the case $d = 2$. 
 Let
$a$ be a smooth
$H^1$ atom on the unit ball, which is smooth and radial, and let $\Omega_N$ be
the lacunary function defined on the unit circle by
$$ \Omega_N(\cos \alpha, \sin \alpha) \equiv G_N(\alpha)=\frac{1}{\sqrt{N}}
\sum_{j=1}^N  e^{2\pi i C^j \alpha},$$
where $C, N$ are large integers.  Roughly speaking, the function $K*a(x)$
has magnitude $\sim N^{-1/2} |x|^{-d}$ whenever $|x| \sim C^j$ for some
$j = 1, \ldots, N$.  This shows that the $L^1$ norm
(and indeed the $L^{1,q}$ quasi-norm for any $q<2$) of $K*a$ grows with $N$,
even though $\Omega$ is in every $L^p$ class, $p < \infty$, uniformly in $N$.
Thus, the best result one can reasonably hope for is that $T$ maps $H^1$
to the Lorentz space $L^{1,2}$, 
 or  the Hardy-Lorentz space $H^{1,2}$, the quasi-norm norm in the latter 
is the $L^{1,2}$ quasinorm of a suitable square-function 
or maximal operator used in the definition of $H^1$  (see \S2 below).

The previous counterexample can be modified to include the 
case $\Omega\in L^\infty$. Take $G_N$ as above, $\eps>0$ and let
$E_{\eps, N}=\{\alpha: |G_N(\alpha)|>{N^\eps}\}$.
Define
$G_{\eps,N}(\alpha)= (G_N(\alpha)(1-\chi_{E_{\eps, N}}(\alpha))$ and
$$\Omega_{\eps,N}(\cos\alpha,\sin\alpha)=G_{\eps,N}(\alpha)-
\frac{1}{2\pi}\int_{-\pi}^{\pi}
G_{\eps,N}(s)ds.$$ Since $G_N$ is in $BMO$ with norm
independent of $N$ we have by the John-Nirenberg inequality that
$|E_{\eps, N}|=O(e^{-cN^\eps})$, for some $c>0$. 
>From this one checks that the $L^1$ norm of $T_{\Omega_N-\Omega_{N,\eps}}a$ 
over the annulus $|x|\sim C^j$ is 
$O(N^{1/2} e^{-cN^\eps}+2^{-j})$, hence negligible.
Since on the other hand $\|\Omega_{N,\eps}\|_\infty\lc N^\eps$ this disproves
a uniform 
$H^1\to L^{1,q}$ estimate  for $q<2/(1+2\eps)$.

\proclaim{Theorem 1.1}
Let $\Omega\in L\log^2 \!L(S^{d-1})$ and assume that
$\int_{S^{d-1}}\Omega d\sigma(\theta)=0$.
Then the operator $T_\Omega$ maps 
$H^1$ to $H^{1,2}$ and also to $L^{1,2}$.
\endproclaim

\remark{Remark 1.2}
In fact we shall see 
 that the $L\log^2 L$ condition can be strengthened to
  an $L\log L$ condition for a Littlewood-Paley square function 
(see Theorem 6.1 below)
\endremark

\medskip

Analogously we may also consider a maximal variant of $T$; 
here no cancellation is imposed. Let $\chi\in C^\infty_0(\Bbb R^d)$
and
$$
\cM_\Omega f(x)=\sup_{h>0} 
\Big|
\int\frac{1}{h^d} \chi(\frac{y}{h})
  \Omega(\frac{y}{|y|}) f(x-y) dy\Big|.
\tag 1.3
$$

\proclaim{Theorem 1.3}
Let $\Omega\in L\log^2\! L(S^{d-1})$. Then $\cM_\Omega$ maps
$H^1$ to $L^{1,2}$.
\endproclaim

Again, a modification of the above example
shows that  $\cM_\Omega$ may fail to map $H^1$ into $L^{1,q}$ for 
$q<2$.

\subheading{Integrals along  curves in the plane}

In this subsection we shall always be working in the plane $\R^2$.  
Let $m > 1$ be a real number; all constants may implicitly depend on $m$.

Define the Hilbert transform $Hf$ and the maximal function $Mf$ 
along the curve $(t,|t|^m)$ by
$$Hf(x)=p.v. \int f(x_1-t,x_2-|t|^m)\ \frac{dt}{t}
\tag 1.4
$$
and
$$
\cM f(x)=
\sup_{h>0}
\Big| \int f(x_1-t,x_2-|t|^m)\frac{1}{h}\eta(\frac t h) dt\Big|;
\tag 1.5
$$
here $\eta$ is a smooth function with compact support.  
These operators are  invariant with respect to the scaling 
$$
(x_1,x_2) \mapsto (t x_1, t^m x_2), \quad t>0.
\tag 1.6
$$

We shall work with the product type  Hardy space on $\Bbb R^2$, considered by
Chang and Fefferman \cite{6} among others; we denote this space by $H^1_\prd$.
Moreover we denote by $H^{1,2}_\prd$ the  product-type  Hardy-Lorentz space 
(see \S2).

\proclaim{Theorem 1.4}
$\cM$  maps $H^1_\prd$ to $L^{1,2}$,
and $H$  maps $H^1_\prd$ to $H^{1,2}_\prd$  and 
 to  $L^{1,2}$. 
\endproclaim

This should be  compared with the
 results of Christ \cite{7} 
 who showed that $M$ and $H$ map the one-parameter Hardy space 
$H^1_{\parab}$ (defined with respect to the dilations (1.6))
to $L^{1,\infty}$, see also  Grafakos \cite{16}. In fact, Christ \cite{7} 
 observes that $H^1_{\parab}$ is not mapped to  $L^{1,q}$ for 
$q<\infty$.

Now let $\gamma=(\gamma_1,\gamma_2)$ be a complex multi-index
with $\Re(\gamma_1), \Re(\gamma_2)\ge 0$, and define the
(pseudo)\-differentiation 
operator $\cD^\gamma$ by
$$ \widehat{\cD^\gamma f} = |\xi^\gamma| \hat f =
|\xi_1|^{\gamma_1} |\xi_2|^{\gamma_2}\hat f.$$
Consider the family of hypersingular operators $H_\gamma$ defined by
$$
H_\gamma f(x_1,x_2)= 
p.v.
\int_{-\infty}^\infty  \cD^\gamma\! f(x_1-t,x_2-|t|^{m})
|t|^{\gamma_1+\gamma_2 m} \frac{dt}{t}.
\tag 1.7
$$
The space  $L^p$  ($1<p<2$)
is not mapped to $L^{p,q}$ if
$q< 2$ (see \cite{8}); moreover this shows that $H$ does not map 
$H^{1}_\prd $ to $L^{1,q}$ 
or any Hardy-Lorentz space  $H^{1,q}$ for any $q<2$. An angular
Littlewood-Paley theory plays a role in this counterexample. 
Grafakos \cite{16}
proved
using the methods in \cite{7}
that for $m=2$, $\gamma_1=0$ and $\Re(\gamma_2)=1-1/p$ the space
$L^p$ is mapped to $L^{p,p'}$
if $1<p\le 2$.
His method surely extends to the general case considered here.

An improved  optimal result is

\proclaim{Theorem 1.5}
Suppose that $\Re(\gamma_1)\ge 0$, $\Re(\gamma_2)\ge 0$ and
$\Re(\gamma_1+\gamma_2)=1-1/p$.
\roster
\itemize{$\bullet$}
 If  $1<p\le 2$ then $H_\gamma$ is bounded from $L^p$ to $L^{p,2}$.

\itemize{$\bullet$}
 If  $p=1$ then $H_\gamma$ is bounded from $H_\prd^1$ to $L^{1,2}$.
\endroster

In both cases the bounds grow at most polynomially in $|\gamma|$.
\endproclaim

\medskip

The following estimate for a localized averaging operator will follow from our proof.
 Let $\eta \in C^\infty_0(\R)$ and define 
$$
A f(x_1,x_2)=\int \eta(t) f(x_1-t, x_2-|t|^m) dt.
\tag 1.8
$$

\proclaim{Corollary 1.6}
Suppose $m\ge 2$. Then $A$ maps $L^{m,2}$ boundedly to the Sobolev space $L^m_{1/m}$.
\endproclaim

\remark{Remarks 1.7}

(i) Suppose that 
$t\mapsto g(t)$ is a smooth curve passing through the origin and suppose that
its curvature  vanishes to at most order $m-2$ at the origin.
  Then the  statement of Corollary 1.8 remains true if $(t,|t|^m)$ is replaced by a $g(t)$ 
provided that $\eta $ is supported in a sufficiently small neighborhood of the origin.

(ii)  In the statements of Theorems 1.4 and 1.5 the curve $(t,|t|^m)$ can be replaced by
$(t,|t|^m \sgn (t))$.

(iii)  A variant of this family $H_\gamma$  was previously  considered
by Stein and Wainger
\cite{27}
in their proof of $L^p$ boundedness  of the Hilbert transform.
They worked
 with
 a  distance function $\rho$, smooth and positive  in $\R^2\setminus\{0\}$
which is  homogeneous of degree 
$1$ with respect to the dilations (1.6)
and considered the analytic family

$$\widetilde H_\alpha f(x_1,x_2)=
\pv \int_{-\infty}^\infty  \rho^\alpha(D) f(x_1-t,x_2-t^{m})
|t|^{\alpha}\frac{dt}{t}.
$$
The result in \cite{27} is that $\widetilde H_\alpha$ is bounded on $L^p$ for 
$\alpha<1-1/p$. Our proof of Theorem 1.3 shows
that this result can be improved to
$\widetilde H_\alpha : L^p\to L^{p,2}$ if $\alpha=1-1/p$, $1<p\le 2$.

(iv) The  principal value singularity $p.v.\, t^{-1} |t|^{\gamma_1+\gamma_2 m}$ in the definition of
$H_\gamma$  can be replaced by 
$\chi^{\gamma_1+\gamma_2 m -1}_+=
\lim_{\eps\to 0} e^{-\eps t}  (\Gamma(\gamma_1+m\gamma_2))^{-1} t_+^{\gamma_1+m\gamma_2-1}$.
This requires only minor changes in the proof of Theorem 1.5.
\endremark

\head{\bf 2. Preliminaries}\endhead
%

\subheading{Notation}
For two quantities $a$ and $b$ we write 
 $a\lc b$ or $b\gc a$ if there exists an absolute 
 positive constant $C$ so that $a\le C b$. 
We shall consistently refer to the homogeneous quasi-norms on Lorentz and
Hardy-Lorentz spaces as ``norms'', even when 
 the triangle inequality with constant $1$ fails.
If $I$ is a (dyadic) cube, then $x_I$ will  denote its center, and
$2^{i_I}$ will denote its side-length.  We somewhat abuse notation and use
$2^s I$ to denote the cube with the same center as $I$ 
and sidelength $2^{s+i_I}$.
The Lebesgue measure of a set $E$ will sometimes be denoted by  $|E|$ and 
sometimes by $\meas(E)$.

\subheading{2A. Hardy spaces}
There are many equivalent characterizations of the  isotropic 
Hardy-spaces  (\cite{13}),
in terms of maximal functions, atomic decompositions and square-functions
(see  \cite{26} for a rather complete treatment).
We shall use  several of them, but most relevant will be the
 characterization via  Littlewood-Paley square-functions, 
 which we choose  as a definition.

Let $\Phi\in \Cal S(\Bbb R^n)$ with the property that
$\widehat \Phi$ is compactly supported  and equal to $1$ in a
 neighborhood of the origin.
Let $\phi_k$ be defined by
$$\widehat \phi_k(\xi)= 
\widehat\Phi(2^{-k-1}\xi)-
\widehat\Phi(2^{-k}\xi)\tag 2.1
$$


Consider the space $\Cal S'_{restr}$ of tempered distributions which
 are {\it  restricted at  $\infty$}; it consists 
of all $f\in \Cal S'$ with the property that 
$f*\phi\in L^r$ for $\phi\in \Cal S$, for sufficiently large $r<\infty$
(we use the terminology of Stein \cite{26, p.123}).
This  choice of the test function space 
allows one to derive 
versions of  the Calder\'on reproducing formula
({\it e.g.} one excludes polynomials which have Fourier transforms  supported at
the  origin). For $0<p,q<\infty$ we
define $H^{p,q}$ as the space  consisting of tempered 
distributions restricted at $\infty$  which satisfy
$$
\|f\|_{H^{p,q}}:=\Big\|\Big(\sum_{k\in\bbZ}|\phi_k*f|^2\Big)^{1/2}
\Big\|_{L^{p,q}}<\infty
\tag 2.2
$$
and write $H^p=H^{p,p}$.  
Using arguments in \cite{13}, \cite{21}
 one can show that the definition does 
not depend on the particular choice  of $\Phi$.
As shown in \cite{21}, \cite{31} some aspects in the classical theory simplify
by assuming (as we do here) that $\widehat \Phi$ has compact support.
In particular for $b>0$, $r>0$  one has the inequality (\cite{21})
$$\sup_{|y|\le 2^{-k}b}|\phi_k*f(x+y)|\le C_{b,r} (M[|\phi_k*f|^r](x))^{1/r}
\tag 2.3
$$
 and (2.3) allows us to take advantage of the 
Fefferman-Stein theorem concerning $L^p(\ell^r)$ 
estimates for the Hardy-Littlewood  maximal function $M$ (\cite{12}).
This carries over to Lorentz-spaces.
Set
$$S_b f(x)=
\Big(\sum_{k\in\bbZ}\sup_{|y|\le b2^{-k}}|\phi_k*f(x+y)|^2\Big)^{1/2}
$$
Since
$\|g\|_{L^{p,q}} \approx \|g^a\|_{L^{p/a,q/a}}^{1/a}$ we obtain that
for $f\in H^{p,q}$
$$
\|f\|_{H^{p,q}}\approx 
\|S_b f\|_{L^{p,q}}.
\tag 2.4
$$
The space  $H^{p,q}$ is  complete quasi-normed space. 
We note that the definition can be extended to 
Hilbert-space valued functions
(in fact when proving estimates we may often reduce to 
finite-dimensional Hilbert spaces with possibly large dimension).

For the purpose of real interpolation consider  the Peetre $K$-functional
$K(t,f,H^{p_0}, H^{p_1})$,
defined for $f\in H^{p_0}+H^{p_1}$  
 as the infimum of
$\|f\|_{H^{p_0}}+t\|f\|_{H^{p_1}}$ over all decompositions
$f=f_0+f_1$ with $f_0\in H^{p_0}$ and $f_1\in H^{p_1}$.
Then a straightforward modification of arguments by Jawerth and Torchinsky
 \cite{19}
yields  the  formula 
$$
K(t,f,H^{p_0}, H^{p_1})\approx 
K(t,S_b f,L^{p_0}, L^{p_1}).
\tag 2.5
$$

Consequently, by (2.4) and (2.5) one identifies $H^{p,q}$ with the real interpolation space
$[H^{p_0}, H^{p_1}]_{\theta, q}$ 
 if $0<\theta<1$ and 
$(1-\theta)/p_0+\theta/p_1=1/p$ (see \cite{2}), and the spaces $H^{p,q}$ can be identified 
with the spaces  in \cite{11}, \cite{15} defined by means of various 
maximal functions or square functions (see \cite{32}).

Let $\{e_k\}$ be an orthonormal basis of $\ell^2$.
>From standard Hardy space theory \cite{26}  we have
$$
\Big \| \sum_k L_k f_k\Big\|_{H^{p,q}} \approx
\Big\| \sum_k \widetilde L_k f_k e_k \Big\|_{L^{p,q}(\ell^2)}
= \Big\|\Big (\sum_k |\widetilde L_k f_k|^2\Big)^{1/2}\Big\|_{L^{p,q}}.
\tag 2.6
$$
where $L_k$, $\widetilde L_k$ denote convolution with $\phi_k$, 
$\widetilde \phi_k$; here  $\widetilde \phi_k$ is as above and
 $\widetilde \phi_k=2^{kd}\widetilde \phi_0(2^k\cdot)$ so that the 
Fourier transform of $\widetilde \phi$ equals one on the support of
 $\widehat\phi$.

Moreover if $E$ is any finite subset of the integers we have
$$
\Big\|\sum_{k\in E} L_k f_k\Big\|_{L^{p,q}}\le C
\Big\|\sum_{k} L_k f_k\Big\|_{H^{p,q}}
\tag 2.7
$$
where $C$ does not depend on $E$. Note, however, that convergence 
in $L^{p,q}$ may not be compatible with convergence in the sense of 
tempered distributions, if $p<1$ or $p=1$, $q>1$.

\subheading {A Littlewood-Paley decomposition}
It is shown in the classical theory that
the above assumptions on $\Phi$ can be substantially weakened.
A general result in this context is in \cite{32}.
To eliminate a number of technical error terms in the proof of 
Theorem 1.1 we shall work with Littlewood-Paley functions 
localized in space, and  in order
to have an analogue of the Calder\'on reproducing formula we will have to
use  a somewhat unusual 
version of the Littlewood-Paley decomposition:

\proclaim{Lemma 2.1}
Let $r$, $N_0$  be  nonnegative integers and let $\eps>0$. 
Then for 
$ s=0,\dots, r$ there are radial  
functions $\Psi_{(s)}$,  
$\psi_{(s)}$ in $C^\infty_0(\Bbb R^d)$ 
with the following properties.

(i) 
$\Psi_s$ is supported on the ball of radius $\eps$ centered at the origin, and
$\widehat{\Psi_s}(\xi)-1=O(|\xi|^{N_0})$ as $\xi\to 0$. Moreover
 $\psi_s=\Psi_s-2^{-d}\Psi_s(2^{-1} \cdot)$ so that the  moments of order 
$\le N_0$
 of $\psi_s$ vanish.

(ii) Define  $\psi^k_s(x)= 2^{kd}\psi_s(2^{k}x)$ and 
let $L^k_s$ be the operator of convolution with $\psi^k_s$. Then for every
tempered distribution $f$ restricted at $\infty$ we have
$$
f = \sum_{k\in \bbZ} L^k_0\cdots L^k_r f;
\tag 2.8
$$
moreover  if $S^0_r$ denotes the operator of convolution with $\Psi_r$ then
$$
f = S^0_r f+ \sum_{k\ge 1} L^k_0\cdots L^k_r f.
\tag 2.9
$$
The convergence in (2.8), (2.9) holds in  the sense of tempered distributions.

\endproclaim

\demo{\bf Proof}
Let $\Psi$ be a radial bump function supported in $\{x:|x|\le 2^{-6r-6}\eps\}$
so that $\widehat \Psi-1 =O(|\xi|^{N+1})$,
 and let $S^k_0$ be the
operator of convolution with $2^{-dk} \Psi(2^{-k} \cdot)$.  Let 
$$
L^k_0 = S^k_0 - S^{k-1}_0.
$$
We recursively define for $s=0,1,\dots,r-1$
$$ 
\align
S^k_{s+1} &= (2\Id - (S^k_s)^2) (S^k_s)^2
\tag 2.10
\\
L^k_{s+1} &= (2\Id - (S^k_s)^2 - (S^{k-1}_s)^2) (S^k_s + S^{k-1}_s)
\tag 2.11
\endalign
$$
and note the identity
$$ S^k_{s+1} - S^{k-1}_{s+1} = (S^k_s - S^{k-1}_s) L^k_{s+1}
\tag 2.12
$$
so that
$S^k_{s+1} - S^{k-1}_{s+1} = L^k_0\cdots L^k_{s+1}$.
One can check inductively that each $S^k_s$ is the  operator  of convolution 
with $2^{kd}\Psi^s(2^k\cdot)$
 where  the radial bump function $\Psi^s$ is supported in 
$\{x:|x|\le 2^{-6(r-s+1)}\eps\}$ and $\widehat{\Psi^s}(\xi)-1=O(|\xi|^{N_0+1})$ as $\xi\to 0$, and that
the operators $L^k_s$, $S^0_s$ have all the desired properties.\qed
\enddemo

\remark{Remark}
We note that 
(2.6) holds if $L_k$, $\widetilde L_k$ are replaced by any of 
 the operators $L_s^k$  above, 
or perhaps by a composition of finitely many such 
operators. This remark holds under the 
condition that the number $N_0$ of vanishing moments is sufficiently large (in dependence of $p$;
specifically we need $N_0\ge n(1/p-1)$).
\endremark

\subheading{Parabolic dilations}
One may define Hardy spaces with respect 
to a nonisotropic dilation structure
\cite{3}. 
 In this paper we need 
to consider such Hardy-spaces on $\bbR^2$ defined with respect to the scaling
$(x_1,x_2) \mapsto (t x_1, t^m x_2)$, for a fixed real number $m>1$.

If we redefine  the function $\phi_k$  to be
$\widehat{ \phi}_k(\xi_1,\xi_2)=
\widehat\Phi(2^{-(k+1)} \xi_1,2^{-(k+1)m}\xi_2)- 
\widehat\Phi(2^{-k} \xi_1,2^{-km}\xi_2)$ 
then the operator of convolution with $\phi_k$ is a Littlewood-Paley projection
to the region 
$|\xi_1| + |\xi_2|^{1/m} \sim 2^k$. We may then define
 $H^p_{\parab}$ as the space of distributions  $f$ 
restricted at $\infty$,
for which
$\| (\sum_k |\widetilde \phi_k*f|^2)^{1/2} \|_p$
is finite.
Similarly one can define parabolic Hardy-Lorentz space and the obvious 
analogues of the statements in the previous subsections remain true.

\subheading{Product type Hardy spaces}
Let $\{L_{k_1,k_2}\}_{k_1,k_2\in \bbZ}$ be a product 
Littlewood-Paley decomposition on $\R^2$, where $L_{k_1,k_2}$ is a 
multiplier with symbol supported in the region 
$\{ (\xi_1,\xi_2): |\xi_1| \sim 2^{k_1}, |\xi_2| \sim 2^{k_2}\}$; 
we may assume  that $L_{k_1,k_2}$ is the operator of convolution with 
$\phi_{k_1}\otimes \phi_{k_2}$ where $\phi_{k_1}$, $\phi_{k_2}$ are as above
(defined on  the real line).

If $0 < p,q < \infty$, 
we define the product Hardy-Lorentz space 
$H^{p,q}_\prd$ to be the quasi-Banach space
which consists of all tempered distributions
restricted at $\infty$ for which
$$ \|f\|_{H^{p,q}_\prd} =
\Big\| \big(\sum_{k_1} \sum_{k_2} |L_{k_1,k_2} f|^2\big)^{1/2} \Big\|_{L^{p,q}}$$
is finite.
We define $H^p_\prd$ to be $H^{p,p}_\prd$.

The formulas for interpolation  of Hardy-Lorentz-spaces remain true; in fact
(2.5) was proved in this context in \cite{19}. Moreover analogues of (2.6), (2.7) remain true for
the operators $L_{k_1,k_2}$. These  can be proved  by using  the theory of product-type singular integral 
operators (see \eg \cite{6}, \cite{14}).

\subheading{2B. Analytic interpolation in Lorentz spaces}
We need a version of a theorem by  Sagher \cite{22}  concerning analytic 
families of operators acting on Lorentz spaces. It has been observed 
in \cite{23} and \cite{16} that Sagher's arguments carry
 over to somewhat more general situations; we now  recall
 the
 version which appeared in \cite{16}.

We denote by $S$  the strip
$S=\{z:0<\Re(z)<1\}$ and by $\overline {S}$ its closure. 
A function $g$ on $\overline{S}$ 
 is said to be of {\it admissible} growth if there is $a<\pi$ so that
$|g(z)|\lc \exp(e^{a|\Im(z)|})$ for $z\in \overline S$.
Let
 $X_0$ and $X_1$ be two Banach spaces, compatible in the sense of interpolation theory, 
and assume that there is a  subspace $W$ of $X_0\cap X_1$ which is dense 
in both $X_0$ and $X_1$. 
For $z\in \overline S$ 
let $\cT_z$ be an operator which maps functions in $W$
to measurable functions on $\Bbb R^n$; $\cT_z$ is then called  an analytic family
if for any $f\in W$ and almost every $x\in \bbR^n$ the function
$z\to \cT_z f(x)$ is analytic in $S$ and continuous and of admissible growth
in   $\overline S$.
Now if 
$$
\|\cT_z f\|_{L^{p_i,q_i}}\le C_i(z)\|f\|_{X_i}, \quad i=0,1,
\tag 2.13
$$
and if $C_i(z)$ is of admissible growth then
the result in \cite{16}  states that 
$T_\theta$ maps the complex interpolation space 
$[X_0,X_1]_\theta$ boundedly to $L^{p_\theta,q_\theta}$;
here
$(1/p_\theta,1/q_\theta)=
(1-\theta)(1/p_0,1/q_0)+\theta(1/p_1,1/q_1)$.
We shall need the following  consequence of this result.
\proclaim{Lemma 2.2}
For $k\in \Bbb Z$ and
$z\in S$ let $T_{k,z}$ be an operator which maps functions 
in $W$
to measurable functions on $\Bbb R^n$ and assume that
$T_{k,z}$ is an analytic family, for each $k$. Suppose that 
for all $f\in W$
$$
\align
\Big \|\sum_{k\in E}|T_{k,i\tau} f|
\Big\|_{L^{1}}
&\le C(i\tau) \|f\|_{X_0}
\tag 2.14
\\
\Big \|\sup_{k\in E}|T_{k,1+i\tau} f|\Big\|_{L^{1,\infty}}
&\le C(1+i\tau) \|f\|_{X_1}
\tag 2.15
\endalign
$$
for any finite subset $E\subset \bbZ$, with admissible constants
$C(i\tau)$, $C(1+i\tau)$. 
Let $0<\theta<1$. Then
$$
\Big \|\Big(\sum_{k\in \Bbb Z}|T_{k,\theta} f|^q\Big)^{1/q}
\Big\|_{L^{1,q}}\lc \|f\|_{[X_0, X_1]_\theta}
\tag 2.16
$$
if $1/q_\theta=1-\theta$.

\endproclaim

\demo{\bf Proof} Fix $\tilde f\in[X_0,X_1]_\theta$ and $E\subset \bbZ$ be finite. 
There are measurable functions
$g_k$ such that
$\sum |g_k(x)|^{q'}\le 1$ and 
$$\Big| \sum_{k\in E}T_{k,\theta} \tilde f(x) g_k(x) \Big|\ge \frac 12
\Big(\sum_{k\in E}|T_{k,\theta} \tilde f(x)|^q\Big)^{1/q}
$$
for almost every $x\in \Bbb R^n$. 
Define 
$g_{k,z}(x)= \frac{g_k(x)}{|g_k(x)|} |g_k(x)|^{q' z}$ if $g_k(x)\neq0$, and
$g_{k,z}(x)=0$ if $g_k(x)=0$.

Now define  an analytic family by $\cT_z f (x)=  \sum_{k\in E}T_{k,z} f(x) g_{k,z}(x) .$
Then the assumptions (2.14-15) imply the boundedness 
of $\cT_{i\tau}$ from $X_0$ to $L^1$ and of $\cT_{1+i\tau}$
 from $X_1$ to $L^{1,\infty}$, with admissible constants. One deduces the boundedness of $\cT_\theta$ from
 $[X_0,X_1]_\theta$ to $L^{1,q}$.
The constants are independent of $E$  and the choice of $\{g_k\}$. 
This implies 
$$
\Big \|\Big(\sum_{k\in E}|T_{k,\theta} \tilde f|^q\Big)^{1/q}
\Big\|_{L^{1,q}}\le C \|\tilde f\|_{[X_0, X_1]_\theta}
$$
with  $C$ being independent of $E$ and $\tilde f$. The finiteness
 assumption on $E$ can 
be removed by  applications of the monotone convergence theorem.
\qed
\enddemo

\subheading{2C. A vector-valued inequality}
We shall use the following observation 
which can  serve as an elementary  substitute 
for the failing $L^p(\ell^1)$ inequality 
for the vector-valued Hardy-Littlewood maximal operator (\cite{12}).
It is just the dual version of a scalar maximal inequality.

\proclaim{Lemma 2.3}
Let $\Phi\in L^1(\Bbb R^d)$ so that for each $\theta\in S^{d-1}$ the function 
$r\mapsto | \Phi(r\theta)|$  is decreasing in $r>0$. 
Let $\{t_k\}_{k\in\Bbb Z}$ be a collection of positive numbers and let
$P_k$ be the operator of convolution with
$t_k^d \Phi(t_k \cdot)$. Then  for $1\le p<\infty$
$$
\Big\|\sum_{k} |P_kf_k|\Big\|_p\le C_p\|\Phi\|_1
 \Big\|\sum_{k}|f_k|\Big\|_p.
\tag 2.17
$$\endproclaim

\demo{\bf Proof}
We may  assume that $\Phi$ is nonnegative.
Then by duality the assertion  follows immediately from the
$L^{p'}$ boundedness of
the maximal operator $w\mapsto \sup_k |P_k w|$; the latter is a 
 consequence of the method of rotation and the bounds for the 
one-dimensional  Hardy-Littlewood operator (see   \cite{26, p.72-73}).\qed
\enddemo

\subheading{2D. Averaging functions in $\boldkey {L}^{\boldkey 1 \boldkey ,
\boldkey q}$}
The triangle inequality fails in $L^{1,q}$ if $q>1$, but the following Lemma, 
proved for $q=\infty$ by 
Stein and N. Weiss \cite{28}, can often   serve as a substitute. For $1<q<\infty$ the statement 
follows from the cases $q=1$ and $q=\infty$ by interpolation.

\proclaim{Lemma 2.4}
Suppose that $\|f_i\|_{L^{1,q}}\le 1$ and $\sum|c_i|\le 1$.
Then $$\Big\|\sum_i c_i f_i\Big\|_{L^{1,q}}\lc
 \sum_i|c_i|(1+\log_+|c_i|)^{1-\frac 1q}.
$$
\endproclaim

\head{\bf 3. A stopping time construction}
\endhead

We shall  use an abstract form of the Calder\'on-Zygmund decomposition, 
in which no nesting or doubling properties are assumed.
The argument is  related to the stopping time construction in
\cite{7}.

\proclaim{Lemma 3.1}
Let $\preceq$, $\seq$ be partial orders on  a set $\Lambda$;  
we also use the notation $\prec$
 synonymously with  $\precneqq$. 
 Let $\Ga$ be a finite subset of $\La$,  let $\nu$ be a 
non-negative measure on $\Ga$, and let
$A:\La\to \Bbb R^+$ be a positive function.

Assume that for each $\gamma\in \Gamma$ and $N>0$ the set
$$\{\la\in \La: A(\la)\le N \text{ and } \gamma\seq \lambda\}
\tag 3.1$$
is finite.

Then one can find a subset $\ffB $ of $\Lambda$ and a map $q:\Ga\to \La$ which have
the following properties.

\roster
\item"{(1)}"  $\ga\seq q(\ga)$ for all $\ga\in \Ga$.

\item"{(2)}"  If  $q(\ga)\notin \ffB$ then $q(\gamma)=\gamma$.

\item"{(3)}" $$\sum_{\la \in \ffB} A(\la) \leq \nu(\Ga)$$

\item"{(4)}"  For all  $\la \in \La$, we have 
$$\nu(\{\ga\in \Gamma: q(\ga) \prec \la, \gamma\seq\la \}) < A(\la).
$$
\endroster
\endproclaim

\demo{\bf Proof} Define
$$
\Lambda_*=  \Gamma\cup \{
\la \in \Lambda: A(\la)\le \nu(\Gamma) \text{ and } \gamma\seq \lambda
\text{ for some } \gamma\in \Gamma\}
\tag 3.2
$$
By the finiteness of $\Gamma$ and the finiteness assumption on the sets (3.1)
the set $\Lambda_*$ is finite.
  Suppose  we have  found $q$ and $\cB$ 
with properties (1)-(4) relatively to $\Lambda_*$ then (1)-(4) are unchanged if
$\Lambda_*$ is enlarged to $\Lambda$.
Hence it suffices to give a proof under the additional assumption that
$\Lambda$ is finite.

We now induct on the cardinality of $\La$.  The lemma is vacuously true when $\La$ is empty, with $\ffB$ being empty and $q$ being the empty function.

Now suppose inductively that $\La$ is non-empty, and that the lemma is true for all sets $\La$ of lesser cardinality.  Choose an element $\la_{max} \in \La$ which is maximal with
respect to the partial ordering $\preceq$, and let 
$\La' = \La - \{\la_{max}\}$.  Define the set $\Ga' \subset \Ga$ by
$$ \Ga' = \Ga \cap \La'$$
if the estimate
$$
\nu(\{\ga\in \Ga: \ga\seq \la_{max}\}) < A(\la_{max})
\tag 3.3
$$
holds, and by
$$ \Ga' = \{ \gamma \in \Ga: \gamma \not \seq \la_{max} \}$$
otherwise.

Now apply the induction hypothesis with $\La$ replaced by $\La'$, $\Ga$ replaced by $\Ga'$, and $A$ and $\nu$ replaced by their restrictions to $\La'$ and $\Ga'$ respectively.  This gives us a set $\ffB' \subset \La'$ and an assignment $q': \Ga' \to \La'$ satisfying analogues 
($1'$)-($4'$) of the desired properties (1)-(4).  

Define the subset  $\ffB$ of $ \La$ by $\ffB = \ffB'$ if (3.3) holds, and $\ffB = \ffB' \cup \{\la_{max}\}$ if (3.3) fails.  Define $q: \Ga \to \La$ by setting $q(\ga) = q'(\ga)$ if $\ga \in \Ga'$, and $q(\ga) = \la_{max}$ if $\ga \in \Ga \backslash \Ga'$.

We now claim that (1)-(4) holds for these choices of $\ffB$ and $q$.  
The claims (1), (2)
are easily verified from ($1'$), ($2'$), 
and the construction of $\ffB$ and $q$.  If (3.3) 
holds then $\ffB = \ffB'$ and (3) follows from ($3'$);
 otherwise, $\ffB = \ffB' \cup \{\la_{max}\}$ and (3) follows from ($3'$), the construction of $\Ga'$, and the failure of (3.3).

It remains to verify (4).  First suppose that $\la \neq \la_{max}$, so that $\la \in \La'$.  
Then (4) follows from ($4'$), because the elements $\gamma$ of $\Ga \backslash \Ga'$ 
satisfy $q(\gamma) = \la_{max}$ and thus cannot contribute to the left-hand side of (4) 
by the maximality of $\la_{max}$.

Now suppose that $\la = \la_{max}$.  If (3.3) holds, then (4) is immediate.  
If (3.3) fails, then by construction the left-hand side of (4) is zero.  
Thus (4) holds in all cases, and the induction step is complete.
\qed
\enddemo

We remark that the finiteness assumption (3.1) may be dropped if one 
is willing to replace the induction by 
transfinite induction ({\it i.e.} use Zorn's lemma). One can then prove this 
lemma for arbitrary $\Lambda$.

\head{\bf 4.
Integrals along  plane curves}\endhead

In this and the next section we shall always be working in the plane
$\R^2$.  We fix a real number $m > 1$, all constants may implicitly
depend on $m$.  We define $H^1_{\parab}$ to be the one-parameter
Hardy space with respect to the scaling (1.6).

The proofs of our results concerning plane curves 
are based 
on the following key estimate.

\proclaim{Proposition 4.1}
For each integer $l$ let  $\eta_l$ be a $C^\infty$ function 
with compact support in $[1/2,2]$ or in $[-2,-1/2]$,
with $C^4$ norms uniformly bounded in $\ell$.

 Let 
 $d\mu_l$ be the  measure defined by
$$ \int f d\mu_l = \int f(x_1-t, x_2 - |t|^m) 2^l \eta_l(2^l t)\ dt.$$
Then for any vector-valued function $F = \{f_{l}\}_{l \in \bbZ}$, 
$$ \Big\| \Big(\sum_l |f_{l} * d\mu_l|^2\Big)^{1/2} \Big\|_{L^{1,2}} \lesssim
\| f \|_{H^1_{\parab}(\ell^2)}.
\tag 4.1$$
We allow the
 $f_l$  themselves  to be Hilbert space valued functions,
and
$|\cdot|$ is then to be interpreted as the Hilbert space norm.

\endproclaim

In the next section, we shall see how this proposition implies $L^{1,2}$ 
and $L^{p,2}$ mapping properties for the Hilbert transform on plane curves 
and similar objects; this will be done by exploiting the fact that the 
$d\mu_l$ have essentially disjoint frequency supports 
if some moment conditions are assumed on the $\eta_l$.  The estimate (4.1) 
should be
 compared with the bound
$$ \big\| \sup_l |f * d\mu_l| \big\|_{L^{1,\infty}} \lesssim
\| f \|_{H^1_{\parab}}$$
proven in Christ \cite{7}. 
Our techniques shall be closely related to those in that paper.

\demo{\bf Proof}
We may decompose $f$ atomically as $f = \sum_I c_I P_I(b_I)$, where the $I$ are
$2^k \times 2^{mk+\vth}$ rectangles with sides parallel to the axes,
and  $k$, $km+\vth$ are integers, $0\le \vth<1$. The $c_I$ 
are non-negative numbers such that $\sum_I c_I \sim 
\| f\|_{H^1_{\parab}(\ell^2)}$, the $b_I$ 
satisfy $\|b_I\|_{L^2(\ell^2)} \lesssim |I|^{-1/2}$, 
and $P_I$ is the projection operator defined by
$$ P_I[b](x) = \Big(b(x) - \frac{1}{|I|} \int_I b(x) dx\Big) \chi_I(x).$$
Note that the definition of
$P_I$ makes sense as acting on  scalar valued functions 
 or on vector-valued functions, as above.
By the   translation trick in \cite{7} (attributed to P. Jones) 
we may assume that the cubes  $I$ 
are dyadic.  Henceforth we shall refer to the $I$ as (parabolic) {\it cubes}.  
It thus suffices 
to show the estimate
$$ \Big\| \Big(\sum_l \Big|\sum_I c_I P_I[b_{I,l}] * d\mu_l\Big|^2\Big)^{1/2} 
\Big\|_{L^{1,2}} \lesssim
(\sum_I c_I) (\sup_I |I|^{1/2})
\Big \| \Big(\sum_l |b_{I,l}|^2\Big)^{1/2}\Big \|_2
\tag 4.2
$$
for arbitrary collections $I$ of cubes, non-negative numbers $c_I$, and 
arbitrary measurable functions $b_{I,l}$.
By limiting arguments it is sufficient to prove the analogue of 
(4.2), where the
sums in $l$ and the sums in $I$ are extended over finite sets
(with bounds independent of the cardinalities). Henceforth we make 
this finiteness assumption.

Fix the   $I$ and $c_I$.
By complex interpolation (Lemma 2.2) it suffices to show that
$$
\Big\| \Big(\sum_l \Big|\sum_I c_I P_I[b_{I,l}] * d\mu_l\Big|^q\Big)^{1/q}
\Big\|_{L^{1,q}} \lesssim
\big(\sum_I c_I\big) \big(\sup_I |I|^{1/q'}\big) \Big\| \Big(\sum_l |b_{I,l}|^q\Big)^{1/q} \Big\|_q
\tag 4.3
$$
holds
for $q = 1$ and $q=\infty$ and all (complex) functions $b_{I,l}$.

When $q=1$, (4.3)
simplifies to
$$ \sum_l \sum_I c_I \big\| P_I[b_{I,l}] * d\mu_l\big\|_1 \lesssim
\big(\sum_I c_I\big) \sup_I \sum_l \|b_{I,l}\|_1,$$
and the claim follows from Young's inequality, the finite mass of $d\mu_l$, and the 
fact that $P_I$ is bounded on $L^1$.  Thus it remains to prove the $q=\infty$ 
endpoint, 
namely
$$
\Big\| \sup_l \big|\sum_I c_I P_I[b_{I,l}] * d\mu_l\big| \Big \|_{L^{1,\infty}}
 \lesssim
\big(\sum_I c_I\big) \sup_I\sup_l |I|\, \| b_{I,l}\|_\infty.
$$

We may assume that
$$
\sup_I\sup_l |I|\, \| b_{I,l}\|_\infty\le 1
\tag 4.4
$$
Writing 
$a_{I,l} = P_I[b_{I,l}]$, we  thus see that $a_{I,l}$ is supported on $I$, has mean zero,
 and $\|a_{I,l}\|_\infty \lesssim |I|^{-1}$, and our task is now to show that
$$
\meas\big(\{ \sup_l |\sum_I c_I a_{I,l} * d\mu_l| \gtrsim \alpha \}\big) \lesssim
\alpha^{-1} \sum_I c_I
\tag 4.5
$$
for all $\alpha > 0$.

Fix $\alpha>0$. We shall
use a sort of 
Calder\'on-Zygmund decomposition 
and will first look at the ``good''
 cubes contributing to a function which is $O(\alpha)$. Let $\cG$ be the family of all
$I$ for which

 $$
M\big(\sum_{I'} c_{I'} \frac{\chi_{I'}}{|I'|}\big)(x) \le  \alpha
\text{ for some $x\in I$};
\tag 4.6
$$
here $M$ is the 
 Hardy-Littlewood maximal operator
with respect to the scaling (1.6).

We consider the contribution of the  cubes in $\cG$  to (4.5).
The $L^\infty$ norm of  $\sum_{I\in \cG} c_I \frac{\chi_I}{|I|}$ is 
$O(\alpha)$, to see this, consider for each  for each $x_0$  
the smallest  cube in $\cG$
 containing $x_0$ and apply (4.6) for this cube.
We  now apply Chebyshev's inequality 
and the standard fact \cite{28} that the  maximal function associated to the curve
 $(t,|t|^m)$ 
is bounded on $L^2$. This yields
$$
\align
\meas\Big(\big\{x:\, &\sup_l\, \big|\sum_{I\in \cG} c_I a_{I,l} * d\mu_l\big|
 \ge \alpha 
\big\}\Big)
\\
&\le \alpha^{-2}\Big\|\sup_l \big|\sum_{I\in\cG} c_I a_{I,l} * d\mu_l\big| \Big\|_2^2
\lc \alpha^{-2}\Big\|\sup_l \sum_{I\in \cG} c_I \frac{\chi_I}{|I|} * 
|d\mu_l|\Big\|_2^2
\\
&\lc \alpha^{-2}\Big\|\sum_{I\in \cG} c_I \frac{\chi_I}{|I|} \Big\|_2^2
\lc \alpha^{-1}\Big\|\sum_{I\in \cG} c_I \frac{\chi_I}{|I|}\Big\|_1
\lc \alpha^{-1} \sum_I c_I.
\tag 4.7
\endalign
$$

Thus we may restrict our attention to the ``bad'' cubes.
By
the Hardy-Littlewood inequality, 
the  $L^{1,\infty}$ norm of
$M(\sum_I c_I\chi_I/|I|)$ is  $O(\sum_I c_I)$, and so by the definition of 
$\cG$
$$ \meas\big(\cup_{I\notin \cG}I\big) \lesssim \alpha^{-1} \sum_I c_I.$$
Let $C>1$ and $CI$
denote the cube expanded by $C$ (with same center as $I$).
 By the Hardy-Littlewood inequality again we  have
$$
\meas\big(\cup_{I\notin \cG} CI\big) \lesssim \alpha^{-1} \sum_I c_I.
\tag 4.8
$$

To complete the proof of (4.5)
we shall prove the stronger square-function estimate
$$
\meas\Big(\big\{x: \big(\sum_l \big|\sum_{I\notin \cG} c_I a_{I,l} * 
d\mu_l(x)\big|^2\big)^{1/2}
 \ge \alpha \big\}\Big)
\lesssim \alpha^{-1}\sum_I c_I.
\tag 4.9
$$

In order to prove (4.9) we use an abstract version of the 
Calder\'on-Zygmund decomposition based on Lemma 3.1.
We first describe the sets $\Lambda$ and $\Gamma$ which occur in this lemma.
If $m\ge 2$ we define  $\Lambda$ as  the set of all dyadic rectangles 
$Q$ of dimensions 
$2^{\sigma} \times 2^{\sigma+(m-1)\tau+\vth}$ 
for integers $\sigma$,  $\tau$ and for  $\vth\in[0,1)$,
 where
$\sigma\le \tau$, and
$(m-1)\tau+\vth$ is the smallest integer $\ge (m-1)\tau$
({\it i.e.} $\vth=0$ if $m$ is an integer).
Note that $\sigma$, $\tau$ and $\vth$ are unique for each 
$Q$ and we shall write 
$\sigma=\sigma(Q)$, $\tau=\tau(Q)$, $\vth=\vth(Q)$. 
If $1<m<2$ we define $\Lambda$ similarly, with the additional 
requirement that  we only admit those $\tau$ for which
 the fractional part of $(m-1)(\tau-\sigma)$ is $<m-1$; this is to ensure that
$\tau(Q)$ is well defined. In both cases the subset $\Gamma$ is  the set of
 parabolic  cubes $I$ for which $c_I\neq 0$ and  which do not belong to
$\cG$; by assumption $\Gamma $ is finite.
  Note that one has $\tau(I) = \sigma(I)$ for parabolic cubes $I$.

We wish to  partially order  the set $\Lambda$  by requiring
$Q\prec Q'$  if  $\tau(Q)< \tau(Q')$; note that then
$Q$ and $Q'$ are incomparable under $\preceq$ if $\tau(Q)=\tau(Q')$ and 
$Q\neq Q'$.
Finally we take set inclusion $\seq$ as the   second partial order  in 
Lemma 3.1.

We define the tendril $T(Q)$ to be the set
$$ T(Q) = \{ x + (t,|t|^m): x \in 2Q, |t| \le 2^{\tau(Q)+2} \}.
\tag 4.10
$$
Note that 
 $ |T(Q)| \sim 2^{\sigma(Q)+m\tau(Q)} + 2^{2\sigma(Q) + (m-1)\tau(Q)}$
for any  rectangle $Q$ parallel  to the axes, and  therefore
$$
|T(Q)| \sim 2^{\sigma(Q)+m\tau(Q)} \quad \text{ for $Q\in \Lambda$. }
\tag 4.11
$$
The function $A(Q)$ in Lemma 3.1 is then defined by
$$A(Q)=\alpha 2^{\sigma(Q) + m\tau(Q)},$$
and the measure $\nu$ is defined by 
$$\nu(\{I\}) = c_I.$$

The finiteness condition in the proof of Lemma 3.1 is easily verified
and we find a map $I\mapsto q(I)$ 
defined  on $\Gamma$
so that
$I\seq q(I)$ and 
$$
\sum\Sb I\in \Gamma\\q(I)\prec Q\\ I\seq Q\endSb c_I<\alpha|T(Q)|
\tag 4.12
$$
for all
$Q\in \Lambda$, and 
$$\meas\big( \cup_{I\in \Gamma} T(q(I))\big)\lc \frac 1{\alpha}\sum_I c_I+
\meas \big(\cup_{I\in \Gamma} T(I)\big);
$$
the latter inequality follows from statements (2) (3), (4) of  Lemma 3.1.
Since for parabolic cubes $I$ the tendril 
$T(I)$ is contained in a fixed dilate of
$I$ and since $\Gamma\cap\cG=\emptyset$ one has  actually
$$\meas\big(\cup_{I\in \Gamma} T(q(I))\big)\lc \frac 1{\alpha}\sum_I c_I,
\tag 4.13
$$
by (4.8).

For any $I$, $l$ we see that $a_{I,l} * d\mu_l$ 
is supported in $T(q(I))$ if
$l < \tau(q(I))$.
In view of (4.13)
the inequality (4.9) 
follows from
$$
\meas\Big(\big \{x: \big(\sum_l |\sum_{I: l \geq \tau(q(I))} c_I a_{I,l} * d\mu_l|^2
\big)^{1/2}\ge \alpha \}\Big) \lesssim
\alpha^{-1} \sum_I c_I.
$$
It suffices by Chebyshev's inequality to prove the $L^2$ estimate
$$\Big\|\Big( \sum_l \Big|\sum\Sb I\in \Gamma\\ l \geq \tau(q(I))\endSb
 c_I a_{I,l} * d\mu_l\Big|^2\Big)^{1/2} \Big\|_2^2 \lesssim \alpha \sum_I c_I.
\tag 4.14
$$

Let $$\Gamma(m)=\{I\in \Gamma:\tau (q(I))=m\}.
\tag 4.15$$
By the triangle inequality it suffices to show
$$\Big\| \Big(\sum_l \big|\sum_{I\in \Gamma(l-s)} c_I a_{I,l} * d\mu_l
\big|^2\Big)^{1/2}
\Big \|_2^2 \lesssim 2^{-s} \alpha \sum_I c_I
$$
for all $s \geq 0$.

Fix $s$.  It then suffices to show that for each $l$
$$ \Big\| \sum_{I\in\Gamma(l-s)} c_I a_{I,l} * d\mu_l \Big\|_2^2 \lesssim
2^{-s} \alpha \sum_{I\in \Gamma(l-s)} c_I
\tag 4.16
$$
for each $l$, since the claim follows by summing in $l$. 
 By scaling (with respect to the parabolic dilations (1.6) and taking into 
account the definition of $\tau(Q)$ we see that it suffices to prove (4.16) 
for $l=0$.  Expanding the left-hand side of (4.16), we reduce to
$$ \sum_{I,I'\in \Gamma(-s)} c_I c_{I'}
|\langle a_{I,0} * d\mu_0, a_{I',0} * d\mu_0 \rangle| \lesssim
2^{-s} \alpha \sum_{I\in\Gamma( -s)} c_I .$$
By symmetry we may assume that $|I'| \leq |I|$.  It then suffices to show that
$$
\sum\Sb I'\in \Gamma(-s)\\ |I'| \leq |I|\endSb c_{I'}
|\langle a_{I,0} * d\mu_0, a_{I',0} * d\mu_0 \rangle| \lesssim
2^{-s} \alpha,
\tag 4.17
$$
for all $I\in \Gamma(-s)$.

Fix $I\in \Gamma(-s)$ with center $x_I$.  $I$ has dimension $2^{\tau(I)} 
\times 2^{m\tau(I)+\vth(I)}$; since $I\seq q(I)$ by Lemma 3.1, (1), 
and $\sigma(q(I))\le \tau(q(I))$
 by definition of $\Lambda$ 
we see  that 
$$
\tau(I) \leq \tau(q(I)) = -s.
\tag 4.18$$

Rewrite the left-hand side of (4.12)
as
$$
\sum_{I': |I'| \leq |I|, \tau(q(I')) = -s} c_{I'}
|\langle a_{I,0} * F, a_{I',0} \rangle|
\tag 4.19
$$
where $F = d\mu_0 * \widetilde {d\mu}_0$
(and  $\widetilde{\,}$  refers to reflection in the argument).
Observe that $F$ is supported on a sector
$$ \{ (x_1, x_2): |x_2| \lc |x_1| \}$$
and obeys the estimates
$$ |\nabla^\alpha F(x)| \lesssim |x|^{-1-|\alpha|}$$
for all multiindices $\alpha$ with $|\alpha| \leq 2$.
>From the size  conditions on $a_{I,0}$, this implies 
$$
 |a_{I,0} * F(x)| \lesssim 2^{-\tau(I)}
$$
and by the moment conditions on $a_{I,0}$
$$ | \nabla^\alpha (a_{I,0} * F)(x)| 
\lesssim 2^{\tau(I)} |x-x_I|^{-2-|\alpha|}, \quad
\text{ if } 
|x-x_I|\ge 2^{\tau(I)+1}, \, |\alpha|\le 1.
$$

This in turn implies from the size and 
moment conditions on $a_{I',0}$ and the assumption $|I'| \leq |I|$ that
$$ | \langle a_{I,0} * F, a_{I',0} \rangle| \lesssim 2^{2\tau(I)}
\diam(I \cup I')^{-3},$$
where the diameter is respect to the Euclidean metric.

Thus it suffices to show that
$$
\sum\Sb I'\in\Gamma(-s)\\ |I'| \leq |I| \endSb
c_{I'} \diam(I \cup I')^{-3}
\lesssim 2^{-2\tau(I)} 2^{-s} \alpha.
\tag 4.20$$

For $\sigma\le -s$, let $\cR_{\sigma,s}$ 
be the set of dyadic rectangles of dimensions
$(2^{\sigma}, 2^{\sigma-(m-1)(s-1)+\vth})$ so that $0\le \vth<1$.
Observe that $\cR_{\sigma,s}$ is a subset 
of $\Lambda$ consisting of rectangles  $R$
with $\tau(R)=-s+1$.
Also 
let $\cW_{\sigma}$ 
be the set of isotropic  dyadic cubes  of dimensions
$(2^{\sigma}, 2^{\sigma})$; then each $W\in\cW_\sigma$ is a union of
$\sim 2^{(m-1)(s-1)}$ rectangles in $\cR_{\sigma,s}$, with disjoint interiors.

If  $I'\in\Gamma(-s)$ with $|I'|\le |I|$ then $I'$ has dimensions
$(2^{\sigma(I')}, 2^{\sigma(I')-(m-1)s})$ and 
$\sigma(I')\le \sigma(I)=\tau(I)$, and therefore 
every such $I'$ is contained in a unique rectangle $R\in\cR_{\tau(I),s}$. 
Since $\tau(q(I'))=-s$ and $\tau(R)=-s+1$ we have
from Lemma 3.1, (4),
$$
\sum\Sb I'\in \Gamma(-s)\\|I'|\le |I|\\ I'\seq R\endSb
c_{I'}\lc
\alpha |T(R)|\lc \alpha 2^{\tau(I)-ms}
$$
and therefore
$$
\align
\sum\Sb I'\in\Gamma(-s)\\ |I'| \leq |I| \endSb
&c_{I'} \diam(I \cup I')^{-3}
\\
&=
\sum_{W\in\cW_{\tau(I)}}\sum\Sb R\in \cR_{\tau(I),s}\\R\subset W\endSb
\sum\Sb I'\in \Gamma(-s)\\|I'|\le |I|\\ I'\subset R\endSb
c_{I'} \diam(I \cup I')^{-3}
\\& 
\lc \alpha 2^{\tau(I)-ms}
\sum_{W\in\cW_{\tau(I)}}
(2^{\tau(I)}+\dist(W,I))^{-3}
\card(\{R\in \cR_{\tau(I),s}:R\subset W\})
\\
&\lc \alpha 2^{\tau(I)-s}
\sum_{W\in\cW_{\tau(I)}}
(2^{\tau(I)}+\dist(W,I))^{-3}
\lc 2^{-2\tau(I)}\alpha 2^{-s}
\endalign
$$
which is (4.20).\qed
\enddemo

\head{\bf 5.  Integrals along plane curves, cont.}\endhead

We now prove Theorems 1.4 and 1.5.
Following \cite{5} we  work with
an angular Littlewood-Paley decomposition.

Let $\zeta\in C^\infty_0(\Bbb R_+)$ so
that $\zeta(s)=1$ if $s\in ((10^m m)^{-1}, 10^m m)$ and define
$Q_l$ by
$$\widehat {Q_l f}(\xi)=
q_l(\xi)\widehat f(\xi)=\zeta(2^{l(m-1)}|\xi_1|/|\xi_2|)\widehat f(\xi).
\tag 5.1
$$
The operators $Q_l$ form a 
Littlewood-Paley family of multipliers supported in
sectors. 
Note that $q_l(\xi)=1$ whenever $\xi$ is 
normal to the curves
$(t, \pm|t|^m)$ if $2^{l-1}\le |t|\le 2^{l+1}$.

Let $\chi_0$ be a smooth  and even function on 
$\Bbb R$ so that $\chi_0(s)=1$ if $|s|\le 1/2$ and
$\chi_0(s)=0$ of $|s|\ge 1$. 
Define $\cP_l$ by
$\widehat {\cP_l f}(\xi)=\chi_0(|(2^{-l}\xi_1, 2^{-lm}\xi_2)|)\widehat f(\xi)$.

Observe that  the multiplier
$q_l$ satisfies the estimates
$\partial^{\alpha} q_l(\xi)=O(|\xi_1|^{-\alpha_1}|\xi_2|^{-\alpha_2})
$
uniformly in $l$.
Therefore by standard product theory we have the estimate

$$\big \| \{(\Id-\cP_l)Q_l f\} \big \|_{H^1_{\prd}(\ell^2)}
\lc
\big \| \{Q_l f\} \big \|_{H^1_{\prd}(\ell^2)}
 \lesssim \| f \|_{H^1_\prd}
\tag 5.2
$$
where
$f$ itself may be a Hilbert-space valued function.

We now consider the maximal function $Mf$.
We show  that
$$ \| \sup_l |d\mu_l*f| \|_{L^{1,2}} \lesssim \|f\|_{H^1_\prd},
\tag 5.3$$
where $d\mu_l$ is a measure as in Proposition 4.1.

Given (5.3) we show the same bound for the nondyadic maximal function by a standard argument.
After a straightforward application of Lemma 2.4 we may assume that
$\eta$ has support in $(-2^{-5}, 2^{-5})$ and vanishes in $(-2^{-6}, 2^{-6})$.
Let $\tilde \eta$ be supported in $\cup \pm(2^{-8}, 2^{-3})$ 
and equal to 
$1$ on  $\cup \pm(-2^{-7}, 2^{-2})$.
We use a Fourier expansion and  write for
$1/2\le s\le 2$
$$\frac 1s\eta(\frac ts)=
\tilde \eta(t)
\sum_{k\in \bbZ} c_k(s)  e^{2 \pi i k t}
$$
where $c_k(s)=O((1+|k|)^{-N})$ uniformly in $s\in [1/2,2]$.
We set $$d\mu_{k,l}= 
 \int f(t, |t|^m) 2^{l} \tilde \eta(2^{l} t) 
  e^{2 \pi i k 2^{l} t} dt.
$$
and $M_kf(x)=\sup_l|f*d\mu_{k,l}|$.
An application of (5.3) shows that $M_k$ maps $H^1$ to $L^{1,2}$ with norm
$O((1+|k|)^4)$ and since $Mf(x)\lc \sum_k(1+|k|)^{-N}M_k f(x)$ we obtain the 
inequality for the nondyadic maximal operator from another 
application of Lemma 2.4.

Now we turn to the proof of (5.3).
As in \cite{5} the idea is to approximate $d\mu_l$ by 
$Q_l(\Id-\cP_l) d\mu_l$ in 
order to apply 
Proposition 4.1
 and (5.2).

Using straightforward integration by parts arguments
we observe that the functions
$\cP_0 d\mu_0$
 and
$ (\Id -\cP_0)(\Id-Q_l) d\mu_0$
are Schwartz functions.  By rescaling this, using (1.6),
we see that the maximal functions
$ \sup_l |f * \cP_l d\mu_l|$ and
$ \sup_l |f * (\Id-\cP_l)(\Id-Q_l)  d\mu_l|$ are dominated by 
nonisotropic version of the  
grand maximal function  (with respect to  (1.6))
 which maps $H^1_{\parab}$ and hence $H^1_\prd$ to $L^1$.  
It thus suffices to show that
$$ \big\| \sup_l |f * (\Id-\cP_l)Q_l d\mu_l|
 \|_{L^{1,2}} \lesssim \|f\|_{H^1_\prd}.
$$
Writing $f_l = (\Id-\cP_l)Q_l f$, we can dominate the left-hand side 
by the $L^{1,2}$ norm of the square-function 
$ (\sum_l |f_l * d\mu_l|^2)^{1/2}$. With this choice of $f_l$
the inequality
$$ \Big\| \big(\sum_l
 |d\mu_l*
f_l|^2\big)^{1/2} \Big\|_{L^{1,2}}\lc\|f\|_{H^1_\prd}
\tag 5.4
$$
follows from  from Proposition 4.1, the embedding 
$ H^1_\prd(\ell^2) \subset  H^1_{\parab}(\ell^2)$ and (5.2).

Now consider the analytic family $H_\gamma$ (and in particular the Hilbert transform
$H=H_0$).
  We may decompose 
$$
H_\gamma f = \sum_l f * d\sigma_l^\gamma
$$
where 
$$
\inn{d\sigma_l^\gamma}{f}=\int f(t,|t|^m) 2^l\chi(2^l t) |t|^{\gamma_1+\gamma_2 m}\frac{dt} t
$$
and $\chi(t)=\chi_0(t)-\chi_0(t/2)$. Note that $\chi$ is an even function.
  The functions 
$\cP_0 d\sigma_0^\gamma$ and
$(\Id-\cP_0)(\Id-Q_l) d\sigma_0^\gamma $ are   Schwartz functions as before, 
but also have mean zero and so their Fourier transforms decay at 0 as well
 as infinity.

 Summing this, we see that $\cD^\gamma \sum_l (\Id-\cP_l)(\Id-Q_l) d\sigma_l$ 
and 
$\cD^\gamma \sum_l \cP_l d\sigma_l$  are 
 standard product type  Calder\'on-Zygmund
kernels
and so convolution with these kernels will preserve $L_p$, $1<p\le 2$  and $H^1_\prd$.  
It thus suffices to show that
$$
 \Big\| \sum_l (\Id-\cP_l)Q_l \cD^\gamma d\sigma_l^\gamma *f \Big\|_{H^{1,2}_\prd} 
\lesssim \|f\|_{H^1_\prd}\quad \text{ if } \Re(\gamma_1+\gamma_2 m)=0
\tag 5.5
$$
and
$$
 \Big\| \sum_l (\Id-\cP_l)Q_l \cD^\gamma d\sigma_l^\gamma *f \Big\|_2
\lesssim \|f\|_2\quad \text{ if } \Re(\gamma_1+\gamma_2 m)=1/2
\tag 5.6
$$
with constants depending polynomially on $\gamma$.

To see (5.6) we note that a standard stationary phase calculation yields
that 
$|\widehat {d\sigma_0^\gamma }(\xi)|\lc (1+|\xi|)^{-1/2}$. 
By scale invariance
we obtain the  uniform  $L^2$
boundedness of the operators with 
convolution kernels
$(\Id-\cP_l)\cD^\gamma d\sigma_l^\gamma$ if $\Re(\gamma_1+m\gamma_2)=1/2$.
The inequality (5.6) 
follows now from the almost orthogonality of the operators $Q_l$.

In order to prove (5.5) 
it suffices to show that
$$ \Big\| \Big(\sum_{k_1,k_2} \big|\sum_l(\Id-\cP_l) Q_l L_{k_1,k_2} f * d\sigma_l^\gamma
 \big|^2\Big)^{1/2} \Big\|_{L^{1,2}} \lesssim
\Big\| \Big(\sum_{k_1,k_2} | L_{k_1,k_2} f |^2\Big)^{1/2}\Big \|_1,
\tag 5.7 $$
by the square function characterization of $H^{1,2}_\prd$; here $L_{k_1,k_2}$ are as in \S2.
For each $k_1$, $k_2$ there are at most $O(1)$ indices $l$ for which 
$(\Id-\cP_l)Q_l L_{k_1,k_2}$ does not vanish, so we may majorize the 
left-hand side of (5.7) by
$$ \Big\| \Big(\sum_{k_1,k_2} \sum_l \big|(\Id-\cP_l)Q_l L_{k_1,k_2} 
f * d\sigma_l^\gamma \big|^2\Big)^{1/2} 
\Big\|_{L^{1,2}}.$$
By Proposition 4.1 we may majorize this in turn by
$$ \big\| \{ 
(\Id-\cP_l)Q_l 
L_{k_1,k_2} f\}_{l,k_1,k_2 \in \bbZ}\big \|_{H^1_{\parab}(\ell^2)}.
$$
But this is bounded by $\|f\|_{H^1_{\prd}}$, 
by standard arguments similar to the proof of (5.2) above.
 This concludes the proof of Theorem 1.5. To see that the 
Hilbert transform $H$ maps $H^1_\prd$
to $L^{1,2}$ we use in addition the product version of inequality (2.7).

Finally we prove Corollary  1.6. 
Define the measures $d\nu_l^\alpha$ by
$$\int{f}{d\nu_l^\alpha}=\int f(t,|t|^m) 
 2^l(\chi(2^l t))\eta(t) 
|t|^{m\alpha}\frac{dt}{t}
$$
and set $d\nu_l=d\nu_l^{1/m}$.
We use duality and prove 
that convolution with   $(\Id-\Delta)^{1/2m} \sum_l d\nu_l$ maps $L^{m'}$ to $L^{m',2}$.

It is easy to see  that for $\theta_1+\theta_2<1$,
$\theta_1\ge 0$, $\theta_2\ge 0$
 the functions
 $ (\Id-\Delta)^{\theta/2}
\sum_l (\Id-\cP_l)(\Id-Q_l) d\nu_l*f$ 
and 
$(\Id-\Delta)^{\theta/2} \sum_l \cP_l d\nu_l*f$  are dominated by
a constant times the  nonisotropic Hardy-Littlewood maximal function of $f$.

Let $\widetilde Q_l=\widetilde q_l(D)$ is defined similarly as $Q_l$ but with
$q_l \widetilde  q_l=q_l$.
 Observe that in view of the compact support of $\eta$ we have 
$d\nu_l^\alpha=0$ if  $l>C_1$ for suitable $C_1$.
Moreover, if $l\le C_1$, we see, using
 the definition of  $Q_l$ and  the Marcinkiewicz 
multiplier theorem that for $\alpha\ge 0$, that
$$\big\|(\Id-\Delta)^{\alpha/2} (\Id-\cP_l) \widetilde Q_l g\big\|_{L^{m',2}}
\lc
\big\|\cD_2^\alpha   Q_l g\big\|_{L^{m',2}}.
$$
Thus it remains to  show that
$$
\big\|\{\cD_2^\alpha Q_l d\nu_l^\alpha*f\}\big\|_{H_\prd^{p,2}(\ell^2)}
\lc \|f\|_{H^{p}_\prd}, \quad \Re(\alpha)=1-1/p,
$$ for $1\le p\le 2$. This is done by a reprise of the arguments above.

\head{\bf 6. Rough homogeneous  kernels: Preliminary reductions}\endhead


Let  $\chi_0$ be  a radial bump function
which is $1$ on $\{x:|x| \leq 1/2\}$ and zero on $\{x: |x| > 1\}$, and
$\chi(x) = \chi_0(x) - \chi_0(x/2)$ is then a function on the unit annulus.
We also denote by $\tilde \chi(t)$ the restriction of $\chi$ 
 to the positive real line $\bbR^+$.

In what follows we shall work with the Littlewood-Paley operators introduced in Lemma 2.1
(with $r=3$) and decompose the identity as
$\Id=\sum_k L^k_0 L^k_1 L^k_2L^k_3 $; we assume that the numbers $N_0$, $\eps$ in Lemma 
2.1 are chosen so that $N_0\ge 100 d$ and $\eps\le 10^{-10d}$.

Let $\del_j$ be the dilation operator
defined by
$$\del_j g(x) = 2^{-jd} g(2^{-j}x),$$
and let $\cA$ be the averaging operator defined by
$$\cA g(x) = C^{-1} \int \tilde \chi(t) t^{-d} g(t^{-1}x) \frac{dt}{t},$$
where $C = \int \tilde \chi(t) \frac{dt}{t}$ is a normalization constant.


Since  $K$ is homogeneous of degree $-d$ we have the decomposition
$$K = \sum_j \del_j \cA[K \chi].
\tag 6.1
$$ 
If the restriction $\Omega$ of $K$  to the unit sphere belongs to
$L \log^2\! L(S^{d-1})$ then 
$K\chi\in L \log^2 \!L(\bbR^d)$ and, 
since standard Calder\'on-Zygmund operators 
 map $L\log^2 L$ to
 $L\log\!L$ 
the  $L\log^2 \!L$  assumption for $K\chi$   is implied by 
$$\Big(\sum_k|L^k_0(K\chi )|^2\Big)^{1/2}\in L\log L.
\tag 6.2
$$

In the present and subsequent section we prove the 
following stronger version of 
Theorem 1.1.

\proclaim{Theorem 6.1}
Let $K$ be homogeneous of degree $-d$ and assume that the restriction $\Omega$
to $S^{d-1}$ is an integrable function satisfying
$\int \Omega d\sigma=0$.
Suppose that (6.2) holds. Then the operator $T_\Omega$ maps 
$H^1$ boundedly to $L^{1,2}$ and also to the Hardy-Lorentz space $H^{1,2}$.
\endproclaim

We also have 
\proclaim{Theorem 6.2}
Let $K_0(r\theta) =\tilde \chi(r)\Omega(\theta) $  and assume 
$\Omega\in L^1(S^{d-1})$ and 
$(\sum_k|L^k_0(K_0)|^2)^{1/2}\in L\log L$. Then $M_\Omega$ maps 
$H^1$ boundedly to $L^{1,2}$.
\endproclaim

We shall prove Theorem 6.1. To prove  Theorem 6.2 we use the argument in \S5 to reduce to 
a  version which involves only dyadic dilations. The proof of the relevant estimate 
for this dyadic maximal operator is similar to the proof of Theorem 6.1 
and therefore omitted.

  Let $\cT$ be the operator
defined by
$$
\cT f =
\sum_j \del_j \cA[K\chi] \,*\,f
\tag 6.3
$$
We now have to show that $\cT$ is bounded from $H^1$ to $H^{1,2}$.
The $H^1\to L^{1,2}$ boundedness follows then from  (2.7) and 
limiting arguments. 
In our proof of (6.3)  we shall  assume that the sum in $j$ is
 actually finite, but prove a bound which is independent of this finiteness assumption. Again a limiting argument proves the general case.

We now decompose 
$f$ in the standard manner
as $f = \sum c_I a_I$, where $c_I$ are nonnegative constants such that
$\sum_I c_I \lesssim \|f\|_{H^1}$, and $a_I$ is an  atom supported
on $I$ with mean zero and
such that $\| a_I\|_\infty \lesssim |I|^{-1}$ 
(\cite{26}). The center of the atom will be denoted by $x_I$ and we may assume that each atom has sidelength $2^{i_I}$ where $i_I$ is an integer.

For technical reasons
we wish to suppress low frequencies in our atoms.
Let
$$\widetilde a_I = \sum_{l \geq -C_0} L^{l-i_I}_0 L^{l-i_I}_1 L^{l-i_I}_2 L^{l-i_I}_3  a_I,
$$
We assume 
$$
\Big\|
\Big(\sum_k|L^k_0(K\chi )|^2\Big)^{1/2}
\Big\|_{L\log L}\le 1
$$
(working with the   norm 
$\|g\|_{L\log ^\gamma\! L}=\inf\{\la>0: \int\frac{|g(t)|}{\la}\log^\gamma 
(e+\frac{|g(x)|}{\la}) dx\le 1\}$) and 
we shall prove that
$$
\Big \| \sum_I
c_I \sum_{j}\del_{j}\cA [
(K \chi)] *  \widetilde a_I \Big\|_{H^{1,2}} \le B \sum_I c_I
\tag 6.4
$$
where $B$ is a constant depending only on $d$.
Now the cancellation of the 
atoms shows that
 $\|a_I - \tilde a_I\|_{H^1} \lesssim 2^{-\eps C_0}$, and so
$$ \big\|f - \sum_I c_I \tilde a_I\big\| \lesssim 2^{-\eps C_0} \|f\|_{H^1}.
\tag 6.5$$
Let $\|\cT\|$ denote the $H^1\to H^{1,2}$ operator-norm, 
which because of our finiteness assumptions is a priori finite.
(6.5) implies
$$
\|\cT f\|_{H^{1,2}}\lc 
2^{-\epsilon C_0}\|\cT\|\|f\|_{H^1} + B \sum c_I.
$$
Therefore, if $C_0$ in the definition of the $\widetilde a_I$ 
is chosen large enough,  this implies that
$\|\cT\|\lc B$.

In what follows we may assume $$\sum c_I\le 1.
\tag 6.6
$$
We now dispose of the contributions when $j \le i_I+ 2C_0$.
We claim this portion is not only in $H^{1,2}$ but is actually in
$H^1$.  Since $H^1$ is a Banach space we may restrict ourselves to a single cube
$I$, so that it suffices to show that
$$ \Big \|\sum_{j\le i_I+ 2C_0} \del_{j}\cA [K \chi] * \tilde a_I 
\Big\|_{H^1}
 \lesssim 1.$$
This we rewrite as
$$ \Big\|\sum_{l \geq -C_0} L^{l-i_I}_0 L^{l-i_I}_1 L^{l-i_I}_2 
[\sum_{j\le i_I+2C_0} 
\del_{j} \cA[K \chi] * L^{l-i_I}_3 a_I] \Big\|_{H^1} \lesssim 1.$$
By the analogue of (2.6) for the Littlewood-Paley operators
$L^k_0 L^k_1 L^k_2$ 
it thus suffices to show
$$ \Big\|\Big(\sum_{l \geq -C_0} \big|\sum_{j\le i_I+ 2C_0} \del_{j} \cA
[K \chi] * L^{l-i_I}_3 a_I\big|^2\Big)^{1/2}\Big \|_1\lc 1.$$
Since the expression inside the norm is supported in 
a fixed dilate of $I$, it suffices by the Cauchy-Schwarz inequality to 
bound 
$$ \Big\|\Big(\sum_{l \geq -C_0} 
\big|\sum_{j\le i_I + 2C_0} \del_{j} \cA[K \chi] * 
L^{l-i_I}_3
a_I\big|^2\Big)^{1/2} \Big\|_2 \lesssim |I|^{-1/2}.$$
By modifying the method of rotations argument
in   \cite{4}
 we see that the operator with convolution kernel
$\sum_{j<i_I+2C_0} \del_{j\le i_I+s} [K \chi]$ is bounded on $L^2$;
hence the above reduces to
$$ \Big(\sum_{l \geq -C_0} \big\|L^{l-i_I}_3 a_I\big\|_2^2\Big)^{1/2}
 \lesssim |I|^{-1/2}.
\tag 6.7$$
But this follows from the $L^2$ estimates on $a_I$ and the almost orthogonality of the $L^{l-i_I}_3$.

We now turn to the contributions  $j> i_I +2C_0$ and we wish to establish
$$
 \Big\|\sum_I c_I\sum_{l \geq -C_0} L^{l-i_I}_0 L^{l-i_I}_1 L^{l-i_I}_2 
[\sum_{j> i_I+2C_0} 
\del_{j} \cA[K \chi] * L^{l-i_I}_3 a_I] \Big\|_{H^{1,2}} \lesssim 1.
$$

We set  $a_{I,l}= L^{l-i_I}_3 a_I$ and let 
 $\{e_j\}$ be  the standard orthonormal basis of unit vectors 
in $\ell^2$.  By the remark following Lemma 2.1 
 we reduce  to show that 
$$
 \Big\|\sum_I c_I\sum_{l \geq -C_0}  L^{l-i_I}_1 
[\sum_{j> i_I+2C_0} 
\del_{j} \cA[K \chi] * L^{l-i_I}_2 a_{I,l}] e_{l-i_I} \Big\|_{L^{1,2}(\ell^2)}
 \lesssim 1.
$$

By Lemma 2.1 we may decompose 
$$K\chi= S^0_1(K\chi)+\sum_{k=1}^\infty L^k_1L^k_0(K\chi).$$
One easily checks that the convolution operator with kernel
$K = \sum_j \del_j \cA[S^0_1 K \chi]$ is a standard 
Calder\'on-Zygmund operator.
Indeed using      the cancellation of the functions $L^{l-i_I}_2 a_{I,l}$
it is easy to see that for a fixed cube $I$
$$
\Big\|\Big(\sum_{l \geq -C_0}\Big|\sum_{j> i_I +2C_0}
\delta_j[\cA S^0_1 (K\chi)] *  L^{l-i_I}_2 a_{I,l}\Big|^2\Big)^{1/2}\Big\|_1\lc 1,
$$
and the resulting $H^1\to L^1(\ell^2)$ inequality follows for this part.

Therefore it suffices to prove that

$$ \Big\| \sum_I c_I \sum_{j > i_I+ 2C_0} \sum_{l \geq -C_0}
L^{l-i_I}_1 \del_{j}\cA(\sum_{k>0}L^k_1 K^k) * L^{l-i_I}_2 a_{I,l}
e_{l-i_I}  
\Big\|_{L^{1,2}(\ell^2)}
\lesssim 1,
\tag 6.8$$
where still  $a_{I,l}= L^{l-i_I}_3 a_I$, and $K^k = L^k_0 (K \chi)$.

We can rewrite the desired estimate for this portion using 
 the identity 
$$L^{m}_1 \del_{j} = \del_{j} L^{j+m}_1.
$$ 
Consequently we have to prove for $q=2$ the inequality
$$
\align
 \Big\| \sum_I \sum_{j> 2C_0+i_I} \sum_{l \geq -C_0}
c_I \del_{j} (L_1^{l-i_I+j} &\cA[\sum_{k>0} L^k_1K^k]) * L^{l-i_I}_2 a_{I,l}e_{l-i_I}
 \Big\|_{L^{1,q}(\ell^q)}\\
&\lc
\sup_I |I|^{1-1/q}  \Big(\sum_l \|a_{I,l}\|_q^q\Big)^{1/q}
\Big\|\big(\sum_k |K^k|^q\big)^{1/q}\Big\|_{L \log^{2-\frac 2q} \!L}
\tag 6.9
\endalign
$$
for {\it arbitrary}  measurable  functions $K^k$ on $\{x:1/4\le |x|\le 4\}$ 
and  $a_{I,l}$ on $CI$. (6.8) follows then by using also  (6.7).

We shall deduce the inequality for $q=2$ from the inequality (6.9)  for $q=1$ 
and the obvious modification of (6.9) for $q=\infty$.

Notice that
$$\align
\big\| L_1^{l-i_I+j} \cA[L^k_1K^k] \big \|_{L^1\to L^1}
&\le \int|\tilde\chi(t)|t^{-d} 
\big\|\psi_1^{l-i_I+j}* t^{-d}\psi^k_1(t^{-1}\cdot)\big\|_1\|K^k\|_1 dt
\\&\lc
2^{-|l-i_I+j-k|}\|K_k\|_1
\tag 6.10
\endalign
$$
where we have  used the cancellation of the 
Littlewood-Paley kernels. The last estimate immediately implies (6.9) 
for $q=1$.
The nontrivial  part concerns the estimate 
for $q=\infty$ which is proved in the next section.
>From these two  estimates  we deduce  (6.9) for $q=2$ 
 by complex interpolation, using 
Lemma 2.2. 
Assuming
$$\Big\|\big(\sum_k |K^k|^2\big)^{1/2}\Big\|_{L \log \!L}\le 1,$$
we consider the analytic family $K_z=\{K^{k}_z\}_{k\in\bbZ}$ defined by
$$
K^{k}_z(x)= K^k(x) |K^k(x)|^{1-2z} \big|K(x)\big|_{\ell^2}^{2z-1}
\big[ \log(e+|K(x)|_{\ell^2})\big]^{1-2z}
$$
if $K^k(x)\neq 0$ and  by
$K^{k}_z(x)=0$ otherwise. Then
$\|K_{i\tau}\|_{L^1(\ell^1)}\lc 1$ and 
$\|K_{1+i\tau}\|_{L\log^2L(\ell^\infty)}\lc 1$. The rest is straightforward.

\head{ \bf 7. Rough homogeneous kernels: The  weak type  estimate}\endhead
We are now proving the analogue of (6.9) for $q=\infty$. In addition to 
(6.6)  we  may also  suppose that
$$
\sup_I\sup_l\|a_{I,l}\|_\infty \le 1,  \qquad
\big\|\sup_k|K_k|\big\|_{L\log^2\!L}\le 1
\tag 7.1
$$
and show that for $\alpha>0$


$$
\meas\Big(\big\{x: \big|\sum_I \sum_{j> 2C_0+i_I} \sum_{l\ge -C_0} c_I
\del_{j}
 (L_1^{l-i_I+j} \cA[\sum_{k>0} L^k_1K^k]) * L^{l-i_I}_2 a_{I,l}e_{l-i_I}
\big|_{\ell^\infty}>\alpha\big\}\Big) \,\lc \, \alpha^{-1}.
\tag 7.2
$$

Let $F=\sum_I c_I \frac{\chi_I}{|I|}$.  Since
$\|F\|_1 \lesssim 1$, we may apply the standard dyadic Calder\'on-Zygmund
decomposition to $F$ at level $\alpha$, and obtain a collection of
disjoint dyadic cubes ${\Cal J} = \{J\}$ 
such that $\sum_J |J| \lesssim \alpha$,
$\int_J F(x)\ dx \lesssim \alpha |J|$,
and such that $F$ is $O(\alpha)$
outside of $\bigcup_J J$.

To every dyadic cube $I$ we assign a nonnegative integer $t_I$ as follows.
If $I$ is not contained in any of the $J$, then $t_I = 0$.  If $I$ is a
subset of a cube $J \in {\Cal J}$, then $t_I$ is chosen so that
the sidelength of $J$ is $2^{t_I}$ times the sidelength of $I$.  One
can view $t_I$ as a stopping time; roughly speaking,
$2^{t_I} I$ is the largest dilate of $I$ on which the mean of $F$ is greater
than $\alpha$, or $I$ if no such dilate exists.

The contribution of the terms in
(7.2) for which $j < i_I+t_I + 2C_0$ is contained inside the exceptional set
$\bigcup_J CJ$, which  has measure $O(\alpha)$.
We can therefore restrict ourselves to the case $j \geq i_I+ t_I + 2C_0$.
We change the summation variable to  $s=j-i_I-t_I \geq 2C_0$.
Thus for the expression
$$
\cE(x)=
\sum_I \sum_{s\geq 2C_0} \sum_l c_I\sum_{k>0}
\del_{i_I+t_I+s}
 (L_1^{l+s+t_I} \cA[L^k_1K^k]) * L^{l-i_I}_2 a_{I,l}(x) e_{l-i_I}
\tag 7.3$$
we have to show that the measure of the set
$\{x:|\cE(x)|_{\ell^\infty}>\alpha\}$ is $O(\alpha^{-1})$.
This  will be estimated by further splitting 
the expression $\cE(x)$ into four pieces  and then  
 by applying of Chebyshev's inequality and
 $L^1$ or $L^2$ estimates for the individual  pieces.

We now describe this splitting.
Let $$M(x)=\sup_{k>0}|K^k(x)|.
\tag 7.4$$

We break up the functions $K^{k}$ into a bounded part and an
integrable part
(this truncation has first been used in \cite{9}).
Let $\eps_0 > 0$ be a constant to be chosen later ($\eps_0=10^{-2}$, say,
 works).
For all  $k$
write
$K^{k} = 2^{\eps_0 (s+l)} K^k_{l,s,I} + R^k_{l,s,I}$, where
$|K^k_{l,s,I}(x)|\le 1$ and 
the remainder $R^k_{l,s,I}$ is the restriction of $K^{k}$ to the set $\{x:M(x)
 \geq 2^{\eps_0 (s+l)}\}$. 
We  split 
$$
\cE(x)= \cE_1(x)+\cE_2(x)+\cE_3(x)+\cE_4(x)
$$
where
$$
\align
\cE_1(x)&= \sum_I \sum_{s\ge 2C_0} \sum_{l\ge -C_0} c_I\sum\Sb
k>0\\ |k-l-s-t_I|\ge s+l\endSb
\del_{i_I+t_I+s}
 (L_1^{l+s+t_I} \cA[L^k_1K^k]) * L^{l-i_I}_2 a_{I,l}(x) e_{l-i_I}
\tag 7.5.1
\\
\cE_2(x)&= \sum_I \sum_{s\ge 2C_0} \sum_{l\ge -C_0} c_I\sum\Sb
k>0\\ |k-l-s-t_I|< s+l
\endSb
\del_{i_I+t_I+s}
 (L_1^{l+s+t_I} \cA[L^k_1R^k_{l,s,I}]) * L^{l-i_I}_2 a_{I,l}(x) e_{l-i_I}
\tag 7.5.2
\\
\cE_3(x)&= \sum_I \sum_{l\ge 2C_0}\sum_{2C_0\le s \le l}
 c_I 
2^{\eps_0(s+l)}
\sum\Sb
k>0\\ |k-l-s-t_I|< s+l
\endSb
\del_{i_I+t_I+s}
 (L_1^{l+s+t_I} \cA[L^k_1 K^k_{l,s,I}]) * L^{l-i_I}_2 a_{I,l}(x) e_{l-i_I}
\tag 7.5.3
\\
\cE_4(x)&= \sum_I \sum_{s\ge 2C_0}\sum_{-C_0\le l< s}
 c_I 
2^{\eps_0(s+l)}
\sum\Sb
k>0\\ |k-l-s-t_I|< s+l
\endSb
\del_{i_I+t_I+s}
 (L_1^{l+s+t_I} \cA[L^k_1 K^k_{l,s,I}]) * L^{l-i_I}_2 a_{I,l}(x) e_{l-i_I}
\tag 7.5.4
\endalign
$$

It suffices to show that for $i=1,2,3,4$ the measure of the set
$\{x:|\cE_i(x)|_{\ell^\infty}>\alpha/4\}$ is $O(\alpha^{-1})$.
By Chebyshev's inequality and the continuous imbedding
$\ell^1\subset\ell^2\subset\ell^\infty$ it suffices to show that
$$
\|\cE_1\|_{L^1(\ell^1)}
+\|\cE_2\|_{L^1(\ell^1)}
+\|\cE_3\|_{L^1(\ell^1)}\lc 1
\tag 7.6
$$
and
$$
\|\cE_4\|_{L^2(\ell^2)} \lc \alpha.
\tag 7.7
$$

The estimation of $\cE_1$ and $\cE_2$ is straightforward.
Since 
$\| (L_1^{l+s+t_I} \cA[L^k_1 K^k]) \|_{L^1\to L^1}\lc 2^{-|k-l-s-t_I|}$ 
we get
$$
\align
\|\cE_1\|_{L^1(\ell^1)}
&\lc 
 \sum_I \sum_{s\ge 2C_0} \sum_{l\ge -C_0} c_I
\sum_{
 |k-l-s-t_I|\ge s+l}
2^{-|k-l-s-t_I|} \| L^{l-i_I}_2 a_{I,l}\|_1
\\
&\lc 
 \sum_I c_I \sum_{s\ge 2C_0} \sum_{l\ge -C_0} 
2^{- s-l} \lc 1.
\tag 7.8 
\endalign
$$

Next, by the definition of $R^k_{l,s,I}$
$$\|L_1^{l+s+t_I} \cA[L^k_1 R^k_{l,s,I}]\|_1\lc
2^{-|k-l-s-t_I|} \int_{x:M(x)\ge 2^{\eps_0(s+l)}}M(x) dx$$
and therefore

$$
\align
\|\cE_2\|_{L^1(\ell^1)}
&\lc 
 \sum_I \sum_{s\ge 2C_0} \sum_{l\ge -C_0} c_I
\sum_{
 |k-l-s-t_I|\le s+l}
2^{-|k-l-s-t_I|} \int_{x:M(x)\ge 2^{\eps_0(s+l)}}M(x) dx
\\
&\lc 
 \sum_I c_I \int |M(x)|\log^2(e+|M(x)|) dx\lc 1.
\tag 7.9
\endalign
$$

The following Lemma is crucial for the estimation of $\cE_3$.

\proclaim{Lemma 7.1}
Suppose that $g$ is a bounded function 
 supported in $\{x: 1/4\le |x|\le 4\}$ and $a$ is supported in a cube $I$ 
with sidelength $2^{i_I}$; moreover  $\|a\|_\infty\le |I|^{-1}$.
Then for $m\ge 0$
$$
\big\|\delta_{i_I+m}[L^{l+m}\cA g]*a\big\|_1\lc 2^{-l/2}\|g\|_\infty
$$
\endproclaim
\demo{\bf Proof}
We may assume $\|g\|_\infty\le 1$.
 Let $\cV_m=\{\nu\}$ be a maximal $2^{-m}$-separated subset of 
unit vectors in $\bbR^d$; its cardinality is $O(2^{m(d-1)})$. We may split
$g=\sum_\nu g_{m,\nu}$ where $g_{m,\nu}$ is supported in the sector
$\{x:|\frac{x}{|x|}-\nu|\lc 2^{-m+10}\}$ 
(and in the annulus where $1/4\le |x|\le 4$).

Now 
$\delta_{i_I+m}[L^{l+m}\cA g]*a$ is supported in a rectangle
of dimensions
$C_1 2^{i_I}\times\dots\times C_12^{i_I}\times C_1 2^{i_I+m}$.
Therefore by the Cauchy-Schwarz inequality
$$
\align
\big\|\delta_{i_I+m}[L_1^{l+m}\cA g]*a\big\|_1
&\lc\sum_{\nu\in \cV_m} 2^{(i_I d+m)/2}
\|\delta_{i_I+m}[L_1^{l+m}\cA g_{m,\nu}]*a\big\|_1
\\
&\lc |I|^{1/2} 2^{md/2}
\Big(\sum_{\nu\in \cV_m} 
\big\|\delta_{i_I+m}[L_1^{l+m}\cA g_{m,\nu}]*a\big\|_2^2\Big)^{1/2}.
\tag 7.10
\endalign
$$
We estimate 
this sum using Plancherel's theorem.
For $\xi\in (\bbR^d)^*$ 
$$
\align
|\widehat {\cA g}_{m,\nu}(-\xi)|
&=\Big|
\int_{r=1/4}^4\int_\theta g_{m,\nu}(r\theta) r^{d-1}
\int\chi(\tau)  e^{i\tau\inn{r\theta}{\xi}} d\tau \, d\theta dr\Big|
\\
&\lc \|g\|_\infty
\int_{1/4}^4\int_{|\theta-\nu|\le 2^{-m+10}} (1+|\inn{\theta}{\xi}|)^{-N} \, d\theta dr.
\\
&\lc 2^{-m(d-1)/2}
\Big(\int_{|\theta-\nu|\le 2^{-m+10}} (1+|\inn{\theta}{\xi}|)^{-2N} \, d\theta
\Big)^{1/2}.
\endalign
$$
Therefore
$$\align
\sum_{\nu\in \cV_m} 
&\big\|\delta_{i_I+m}[L_1^{l+m}\cA g_{m,\nu}]*a\big\|_2^2
\\
&\lc 
2^{-m(d-1)}  \sum_{\nu\in \cV_m} 
\int\big|\widehat {\psi_1^{l+m}}(2^{i_I+m}\xi)\big|^2
\int_{|\theta-\nu|\le 2^{-m+10}} 
(1+|\inn{\theta}{2^{i_I+m}\xi}|)^{-2N} d\theta\,
|\widehat a(\xi)|^2 d\xi
\\
&\lc 2^{-m(d-1)} 
\int\big|\widehat {\psi_1^{l+m}}(2^{i_I+m}\xi)\big|^2
\int_{S^{d-1}}
\big| (1+|\inn{\theta}{2^{i_I+m}\xi}|)^{-2N}  d\theta\,
|\widehat a(\xi)|^2 d\xi
\\
&\lc 2^{-m(d-1)} 
\int\big|\widehat {\psi_1}(\frac{\xi}{2^{l-i_I}})\big|^2
\,\min\{1, 2^{-i_I-m}|\xi^{-1}|\}\,
|\widehat a(\xi)|^2 d\xi
\\
&\lc 2^{-m(d-1)} 2^{-(m+l)} \|\widehat a\|_2^2 \lc
 2^{-md-l}|I|^{-1},
\tag 7.11
\endalign
$$
by Plancherel's theorem and the estimate 
$|\widehat {\psi_1}(\xi)|\lc \min\{|\xi|^2,
|\xi|^{-2}\}$.

The asserted estimate follows from (7.10) and (7.11).\qed
\enddemo

We now estimate the $L^1(\ell^1)$ norm of $\cE_3$.
To apply Lemma 7.1 we note that
$L^{l-i_I}_2 a_{I,l}$ is supported in a fixed dilate of $I$ and 
$\|L^{l-i_I}_2 a_{I,l}\|_\infty \lc |I|^{-1}$. 
Moreover 
$\|L^k_1K^k_{l,s,I}\|_\infty\lc 1$, uniformly in $k,l,s,I$.
Hence  
$$
\align
\|\cE_3\|_{L^1(\ell_1)}\lc
 \sum_I c_I \sum_{l\ge 2C_0}\sum_{2C_0\ge s\le l}
2^{\eps_0(s+l)} 
\sum\Sb
k>0\\ |k-l-s-t_I|< s+l
\endSb
2^{-l/2} \|L^k_1K^k_{l,s,I}\|_\infty
\lc 1.
\tag 7.12
\endalign
$$

Finally we turn to the estimation of $\|\cE_4\|_{L^2(\ell^2)}$.
We first observe the basic estimate
\proclaim{Lemma 7.2}
$$
\Big\|\sum_I c_I \frac{\chi_{2^{t_I} I}}{|2^{t_I} I|}\Big\|_2 \lesssim \alpha^{1/2}.$$
\endproclaim

\demo{\bf Proof}
  Consider first those cubes $I$ for which $t_I = 0$.  
It is easy to see that this contribution is bounded pointwise by 
$\min(F, C\alpha)$ for some constant $C$, 
and so the claim follows since $\|F\|_1 \lesssim 1$.

Now consider the cubes $I$ for which $t_I > 0$.  
This part is majorized pointwise by
$$\Big\|\sum_J \chi_{CJ}\Big\|_2 \lesssim \Big\|\sum_J \chi_J\Big\|_2 
= \Big(\sum_J |J|\Big)^{1/2} \lesssim \alpha^{1/2},$$
where for the first inequality we have used 
Lemma 2.3. \qed
\enddemo

The claimed estimate for $\cE_4$ will follow from

\proclaim{Lemma 7.3}
Let $g_I$ be bounded and  supported on $\{x:1/4\le |x|\le 4\}$ and set
$b_{I,l}= L_2^{l-i_I} a_{I,l}$. Assume $l\ge -C_0$, $s\ge 0$.
Then for suitable $\eps >0$
$$
\Big\|\sum_I c_I
\del_{i_I+t_I+s}(L_1^{l+s+t_I} \cA g_I)* b_{I,l}
\Big\|_2\lc  \sup_I\|g_I\|_\infty 2^{-s\eps}\alpha^{1/2}.
$$
\endproclaim

\demo{\bf Proof}
This inequality is closely related  to one in \cite{25} and we shall adapt the proof here.
Let $\cV_s=\{\nu\}$ be a maximal $2^{-s}$-separated subset of the
unit sphere $S^{d-1}$;  the cardinality of this set is
$O(2^{(d-1)s})$.  We decompose $g_I = \sum_\nu g_{I,\nu}$, 
where each $g_{I,\nu}$ is a bounded function on the sector
$$
\fS_\nu^s=
\{x : 1/4\le |x|\le 4, 
\angle( x,\nu) \le 2^{-s} \};
\tag 7.13
$$
here we used  $\angle(x,\nu)$ to denote the angle
$x$ and $\nu$ make at the origin.

We introduce a localization in Fourier space to a conic
neigborhood of the hyperplane perpendicular to $\nu$, namely
$$\Sigma_\nu^s=\{\xi: |\inn \xi \nu|\le 2^{- s/2}|\xi|\}
$$
 (The exact choice of aperture $2^{-s/2}$ is unimportant as long as 
 it is well between $2^{-s}$ and 1).
We define the
multiplier $Q_\nu^s$ whose symbol $m_\nu$ is homogeneous of degree 0, 
and equals 1 on $\Sigma_\nu^s$ and vanishes outside  a slight widening of $\Sigma_\nu^s$.

We then reduce to showing that

$$
\Big\|\sum_I c_I \sum_{\nu\in \cV_s}
Q_\nu^s \del_{i_I+t_I+s}(L_1^{l+s+t_I}  \cA g_{I,\nu})* b_{I,l}
\Big\|_2\lc \sup_I\|g_I\|_\infty 2^{-s\eps}\alpha^{1/2}
\tag 7.14
$$
and, for fixed $\nu$,
$$
\Big\|\sum_I c_I
(\Id-Q_\nu^s)
\del_{i_I+t_I+s}(L_1^{l+s+t_I} \cA g_{I,\nu})* b_{I,l}
\Big\|_2\lc \sup_{I,\nu}\|g_{I,\nu}\|_\infty 2^{-s N}\alpha^{1/2}
\tag 7.15
$$
where $N\le N_0/10$ (recall that $N_0\ge 100 d$).
The estimate (7.15) is favorable if $N>d-1$.

To prove (7.15) we show the estimate
$$
\big|(\Id-Q_\nu^s)
\del_{i_I+t_I+s}(L_1^{l+s+t_I} \cA g_{I,\nu})(x)\big|\lc 
\|g_{I,\nu}\|_\infty
 2^{-sN}\frac{2^{-(i_I+t_I)d}}
{(1+2^{-(i_I+t_I)}|x|)^{N}}
\tag 7.16
$$
for all $\nu\in \cV_s$.
>From (7.16) we may estimate
$$
\big|(\Id-Q_\nu^s)
\del_{i_I+t_I+s}(L^{l+t_I+s} \cA g_{I,\nu})* b_{I,l}\big|\lc 2^{-Ns} 
H_{I}*  \frac{\chi_{2^{t_I} I}}{|2^{t_I} I|}
$$
where  $H_{I}$ is the $L^1$ dilate of a radially  decreasing  $L^1$ function.
By Lemma 2.3 and Lemma 7.2 
the left hand side of (7.15) is dominated by
$$
2^{-sN}
\Big \|\sum_I c_I \frac{\chi_{2^{t_I} I}}{|2^{t_I} I|}\Big\|_2\lc 
2^{-sN}\alpha^{1/2}.
$$

We now show  (7.16).
 Fix $\nu$.  Rescaling so that $i_I+t_I+s = 0$, it suffices to show that
$$ | L^j_1 (\Id - Q^s_\nu) \cA h(x) | \lesssim 2^{-(N+d)s}
\|h\|_{L^\infty(\fS_\nu^s)} (1 + |x|)^{-N}$$
for all $j \ge l+t_I+s\ge s$ and all bounded $h$ supported on $\fS_\nu^s$.

Fix $j, x$.  We expand the left-hand side as
$$ \Big|(2\pi)^{-d} \int_{\fS^\nu_s} h(z) \iiint
(1-m_\nu(\xi)) e^{ i \langle \xi, x-2^{-j} y-tz\rangle}
\psi_1(y) \tilde \chi(t)\, d\xi dy
\frac{ dt}t\,  dz \Big|$$
where the moments of $\psi_1$ vanish up to order $N_0$ and $\tilde \chi$ is supported where 
$1/4\le t\le 4$.
The decay in $x$ follows from the fact that the phase is 
non-stationary in the $\xi$ variable when $|x| \gg 1$.

  Now we demonstrate 
the $2^{-Ns}$ bound; we may assume that $|x|\ll 2^{s/5}$.  Since $h$ is supported in
$\fS_\nu^s$ and  $m_\nu$ equals $1$ on $\Sigma_\nu$  we see that 
for each $|\xi| \gtrsim 2^j$, the phase is non-stationary in the
  $t$ variable (with a gradient of at least $2^{\eps s}$).  
For $|\xi| \lesssim 2^j$ one picks up a loss of $(2^j/|\xi|)^C$, 
but this can be compensated for by the moment conditions on $\psi_1$, since $j\ge s$.

To show 
(7.14)
we use the fact that the $Q^s_\nu$ have some weak
orthogonality.  More precisely, we have for any functions $f_\nu$ that
$$\Big \|\sum_\nu Q^s_\nu f_\nu\Big\|_2^2 \lesssim 2^{-\eps s} 2^{(d-1)s}
\sum_\nu \|f_\nu\|_2^2;
\tag 7.17$$
as in \cite{25} 
this estimate is easily proven from Plancherel's theorem, the Cauchy-Schwarz
inequality, and geometrical considerations.  Because of this orthogonality,
and Lemma 7.2,
it now suffices to show that
$$ \Big\| \sum_I c_I \del_{i_I+t_I+s} \cA g_{I,\nu} * a_I \Big\|_2
\lesssim 2^{-(d-1)s}
\Big\| \sum_I c_I \frac{\chi_{2^{t_I} I}}{|2^{t_I} I|}\Big\|_2,
\tag 7.18$$
uniformly in $\nu\in \cV_s$.

Fix $\nu$.  
Let $R_\nu^s$ be the  rectangle centered at the origin, with 
dimensions $C_1 2^{-s}\times\dots\times C_1 2^{-s}\times C_1$ so that the long side is 
parallel to $\nu$. Then,  if $C_1$ is chosen large enough
there is the uniform pointwise estimate
$$
\big|\del_{i+t_I+s} [\cA g_{I,\nu}] * a_I\big | \lesssim 2^{-s(d-1)}\|g_{I,\nu}\|_\infty
\, \delta_{i_I+t_I+ s}\big(\frac {\chi_{R_\nu^s}}{|R_\nu^s|}\big) *
\frac {\chi_{2^{t_I}I}}{|2^{t_I}I|}.
$$ 
Thus (7.18) follows from Lemma 2.3. This completes the proof of (7.14) and the Lemma.
\qed
\enddemo

The estimate (7.7) is an immediate consequence of Lemma 7.3. The estimate
(7.6) holds by 
 (7.8), (7.9) and (7.12) and thus we have proved
the asserted weak type inequality.

\enddocument